%% file: a12.tex
\documentclass{amsart}
\usepackage{amssymb,xspace}
\usepackage{graphics,graphicx,color}
\usepackage[small]{diagrams}
\input{command}\theorems
\usepackage{hyperref}

\graphicspath{{images/}}
\makeatletter
\def\input@path{{images/}}
\makeatother

\DMO{\inn}{inn}
\def\fl{\mathrm{fl}}
\DMO{\temp}{mod}\let\mod\temp 
\def\a{A} 
\def\b{B} 

\begin{document}

\input{title.tex}\maketitle
\tableofcontents
\input{introduction.tex}
\input{s01.tex} 
\input{s02.tex} 
\input{s03.tex} 
\input{s04.tex} 
\input{example.tex}

\bibliography{../tex/fullbib}
\bibliographystyle{../tex/hamsplain}
\end{document}

%% file: command.tex
\def\mat#1{\ensuremath{#1}\xspace}
\def\makedef#1#2{\expandafter\gdef\csname #1\endcsname{#2}}
\def\makelet#1#2{\expandafter\let\csname #1\expandafter\endcsname\csname #2\endcsname}
\def\defmath#1{\makelet{temp@#1}{#1}\makedef{#1}{\mat{\csname temp@#1\endcsname}}}

\def\defbb#1{\makedef{c#1}{\mat{\mathbb{#1}}}}		
\def\defcal#1{\makedef{l#1}{\mat{\mathcal{#1}}}}	
\def\bbcal#1{\defbb{#1}\defcal{#1}}
\bbcal A
\bbcal B
\bbcal C
\bbcal D
\bbcal E
\bbcal F
\bbcal G
\bbcal H
\bbcal I
\bbcal J
\bbcal K
\bbcal L
\bbcal M
\bbcal N
\bbcal O
\bbcal P
\bbcal Q
\bbcal R
\bbcal S
\bbcal T
\bbcal U
\bbcal V
\bbcal W
\bbcal X
\bbcal Y
\bbcal Z

\def\al{\mat{\alpha}}
\def\be{\mat{\beta}}
\def\ga{\mat{\gamma}}
\def\Ga{\mat{\Gamma}}

\def\eps{\mat{\varepsilon}}

\def\hi{\mat{\chi}}

\def\La{\mat{\Lambda}}
\def\la{\mat{\lambda}}

\def\Om{\mat{\Omega}}
\def\om{\mat{\omega}}

\def\si{\mat{\sigma}}

\def\te{\mat{\theta}}

\defmath{eta}
\defmath{xi}
\defmath{mu}
\defmath{nu}
\defmath{Phi}
\defmath{Psi}
\defmath{pi}
\defmath{rho}

\def\mrm@#1{\mat{\mathrm{#1}}}

\def\deffrak#1{\makedef{g#1}{\mat{\mathfrak{#1}}}} 
\deffrak b
\deffrak c
\deffrak d
\deffrak D

\deffrak h
\deffrak m
\deffrak n
\deffrak p
\deffrak P
\deffrak q
\deffrak r
\deffrak s
\deffrak t
\deffrak u
\deffrak v

\def\DMO{\DeclareMathOperator}

\DMO{\Hom}{Hom}
\DMO{\RHom}{RHom}
\DMO{\lHom}{\lH\mathit{om}}
\DMO{\Ext}{Ext}
\DMO{\lExt}{\lE\mathit{xt}}
\DMO{\End}{End}
\DMO{\Aut}{Aut}
\DMO{\Fun}{Fun}
\DMO{\Tor}{Tor}
\DMO{\ext}{ext}

\DMO{\Ob}{Ob}
\DMO{\Mor}{Mor}
\DMO{\im}{im}
\DMO{\coim}{coim}
\DMO{\coker}{coker}
\DMO{\Arr}{Arr}

\DMO{\Id}{Id}
\DMO{\id}{id}
\DMO{\add}{add} 
\DMO{\ind}{ind} 
\DMO{\pro}{pro} 
\DMO{\Map}{Map} %
\DMO{\Iso}{Iso} %
\DMO{\Isom}{Isom}%
\DMO{\Ind}{Ind}

\DMO{\Presh}{Presh}

\DMO\coalg{Coalg}
\DMO{\Rep}{Rep}
\DMO{\Cor}{Cor}

\DMO{\Mod}{Mod}
\DMO{\rad}{rad}
\DMO{\soc}{soc}
\DMO{\ann}{ann}

\DMO{\Spec}{Spec}
\DMO{\Specm}{Specm}
\DMO{\spec}{Spec}
\DMO{\Proj}{Proj}
\DMO{\supp}{supp}
\DMO{\coh}{coh}
\DMO{\Qcoh}{QCoh}
\DMO{\QCoh}{QCoh}
\DMO{\Pic}{Pic}
\DMO{\Div}{Div}
\DMO{\ch}{ch}
\DMO{\Hilb}{Hilb}
\DMO{\Fitt}{Fitt}
\DMO{\Quot}{Quot}

\DMO{\Gras}{Gr}
\DMO{\Flag}{Flag}

\DMO{\cone}{cone}
\DMO{\Tw}{Tw}
\DMO{\rank}{rk}
\DMO{\rk}{rk}
\DMO{\codim}{codim}
\DMO{\cov}{cov}
\DMO{\sgn}{sgn}
\DMO{\td}{td}
\DMO{\GL}{GL}
\DMO{\SL}{SL}

\DMO\Der{Der}
\DMO\der{Der}
\DMO\coder{Coder}
\DMO{\diag}{diag}
\DMO{\HMod}{HMod} 
\DMO{\ad}{ad}
\DMO{\Ad}{Ad}
\DMO*{\colim}{colim}
\DMO*{\hocolim}{hocolim}
\DMO*{\holim}{holim}
\DMO{\Ho}{Ho}

\DMO{\har}{char}
\DMO{\sk}{sk}
\DMO{\cosk}{cosk}
\DMO{\Gal}{Gal}
\DMO{\tr}{tr}
\DMO{\Tr}{Tr}
\DMO{\Sh}{Sh}
\DMO{\Is}{Is} 
\DMO{\Hol}{Hol} 
\DMO{\Lie}{Lie} 
\DMO{\Res}{Res} 
\DMO{\irr}{irr} %
\DMO{\Irr}{Irr} %
\DMO{\Exp}{Exp} %
\DMO{\Log}{Log} %
\DMO{\Pow}{Pow}
\DMO{\pow}{pow}
\DMO{\mult}{mult} %
\DMO{\height}{ht} %
\DMO{\wt}{wt}
\DMO{\Vect}{Vect}
\DMO{\hd}{hd} 

\newcommand{\GIT}{/\!\!/}

\def\dd{\mat{\partial}}
\def\iso{\simeq}
\def\ts{\otimes}
\def\tl#1{\mat{\tilde{#1}}}

\def\sb{\subset}

\def\ms{\backslash} 
\DMO{\Eig}{Eig} 

\def\op{^{\rm op}}
\def\inv{^{-1}}
\def\dual{^\vee}

\def\angs#1#2{\mat{\left\langle #1 \mid #2\right\rangle}}
\def\set#1{\mat{\{ #1\}}}
\def\sets#1#2{\mat{\{ #1 \mid #2\}}}

\def\emb{\hookrightarrow}
\def\mto{\mapsto}
\def\arr{\futurelet\test\arrtest}
\def\arrtest{\ifx^\test\let\next\arra\else\let\next\arrb\fi\next}
\def\arra^#1{\xrightarrow{#1}} \def\arrb{\to}

\def\arrowsD{
\def\mto{{\:\vrule height .9ex depth -.2ex width .04em\!\!\!\;\ar}}
\def\ar{{\:\vrule depth -.52ex height .60ex width 0.85em\;\!\!\rhla\,}}
\def\arr{\futurelet\test\arrtest}
\def\arrtest{\ifx^\test\let\next\arra\else\let\next\arrb\fi\next}
\def\arra^##1{\rTo^{##1}} \def\arrb{\ar}
\def\emb{\futurelet\test\embtest}
\def\embtest{\ifx^\test\let\next\emba\else\let\next\embb\fi\next}
\def\emba^##1{\rInto^{##1}} \def\embb{{\:\rthooka\!\!\!\ar}}
\newarrow{Eq}=====
\def\rrarr{\pile{\rTo\\ \rTo}}
\def\lrarr{\pile{\rTo\\ \lTo}}   
\newarrow{ShortTo}{}{}-->
}

\def\arrowsDStandard{
\newarrow{TeXto}----{->}
\newarrow{TeXinto}C---{->}
\newarrow{TeXonto}----{->>}
\newarrow{TeXdashto}{}{dash}{}{dash}{->}
\newarrow{Eq}=====
\def\ar{\rightarrow}
\def\emb{\futurelet\test\embtest}
\def\embtest{\ifx^\test\let\next\emba\else\let\next\embb\fi\next}
\def\emba^##1{\rTeXinto^{##1}} \def\embb{\hookrightarrow}
\def\rrarr{\pile{\rTo\\ \rTo}}
\def\lrarr{\pile{\rTo\\ \lTo}}   
}

\newif\ifukr\ukrfalse
\newif\ifrus\rusfalse
\newif\ifger\gerfalse

\def\theorems{
\newcounter{nthr} 
\numberwithin{nthr}{section}
\let\theHnthr\thenthr
\newtheorem{thr}[nthr]{Theorem}
\newtheorem{prp}[nthr]{Proposition}
\newtheorem{lmm}[nthr]{Lemma}
\newtheorem{crl}[nthr]{Corollary}
\newtheorem{clm}[nthr]{Claim}
\newtheorem{conj}[nthr]{Conjecture}
\theoremstyle{definition}
\newtheorem{dfn}[nthr]{Definition}
\newtheorem{rmr}[nthr]{Remark}
\newtheorem{exm}[nthr]{Example}
\newtheorem{claim}[nthr]{Claim}}

\def\ub#1{\mat{\overline{#1}}}  

\def\eprint#1#2{
\expandafter\ifx\csname eprint@#1\endcsname\relax#1:#2\else%
\expandafter\def\csname\csname eprint@#1@idname\endcsname\endcsname{#2}%
\csname eprint@#1\endcsname\fi}

\def\defArchive#1#2#3{ 
\makedef{eprint@#1}{#3}\makedef{eprint@#1@idname}{#2}}

\defArchive{arxiv}{id}{\href{http://arXiv.org/abs/\id}{{\tt arXiv:\id}}}
\defArchive{numdam}{id}{\href{http://www.numdam.org/item?id=\id}{Numdam:\id}}


%% file: title.tex
\title[Crepant resolutions and brane tilings]
{Crepant resolutions and brane tilings I: Toric realization}%

\author{Sergey Mozgovoy}%


\email{mozgov@math.uni-wuppertal.de}%


\begin{abstract}
Given a brane tiling, that is, a bipartite graph on a torus,
we can associate with it a singular $3$-Calabi-Yau variety.
In this paper we study its commutative and non-commutative
crepant resolutions. We give an explicit
toric description of all its commutative crepant resolutions.
We also explain how the McKay correspondence in dimension $3$
can be interpreted using brane tilings.
\end{abstract}

%% file: introduction.tex
\section{Introduction}
The main goal of this paper is to give an explicit
toric description of all possible crepant resolutions of singular
$3$-Calabi-Yau varieties arising from brane tilings (see Section \ref{pre:brane}).
All these crepant resolutions can be constructed
as moduli spaces of representations of some
quiver with relations \cite[Theorem~15.1]{Ishii2}. 
It follows from the construction of these
moduli spaces that they are toric varieties.
We give an explicit description of the corresponding
fan. 

With any brane tiling we can associate a quiver potential $(Q,W)$
(see Section~\ref{pre:brane}).
It turns out that under
certain consistency conditions on the brane tiling
the corresponding quiver potential algebra $\cC Q/(\dd W)$
is a $3$-Calabi-Yau algebra \cite{MR2,Broomhead1}.
The singular Calabi-Yau variety mentioned above
is isomorphic to the spectrum of the center of $\cC Q/(\dd W)$.
We will show that this variety is a normal Gorenstein toric
variety and that $\cC Q/(\dd W)$ is its non-commutative
crepant resolution \cite{Bergh1}.
Related questions are studied in \cite{Bock2}.

Using this fact, we can
apply the result of Van den Bergh \cite[Theorem~6.3.1]{Bergh1}, which says
that certain moduli spaces of representations of a non-commutative
crepant resolution give rise to (commutative) crepant
resolutions and, moreover, the derived categories of
commutative and non-commutative crepant resolutions
are equivalent. As Van den Bergh mentions, his result
is a generalization of the well-known approach of
\cite{BKR} to the McKay correspondence.
In the case of brane tilings, the moduli spaces are
$\lM_\te=\lM_\te(\cC Q/(\dd W),\al)$ -- the moduli
spaces of $\te$-semistable $\cC Q/(\dd W)$-representations
of dimension $\al=(1,\dots,1)\in\cZ^{Q_0}$, where
$\te\in\cZ^{Q_0}$ is \al-generic.
A direct proof of the smoothness of these moduli spaces 
was given by Ishii and Ueda \cite{Ishii1}.
They also proved the derived equivalence \cite{Ishii2}
by using some tricky modifications of brane tilings.

The moduli space $\lM_\te$ has a natural action of a certain
$3$-dimensional torus.
We will prove that every orbit in $\lM_\te$ is
determined by its cosupport -- the subset of arrows
in $Q_1$ inducing zero action on the representations
from the orbit.
It was proved in \cite[Lemma~6.1]{Ishii1} that the cosupports
of $2$-dimensional orbits are perfect matchings.
We can give similar descriptions for the cosupports
of $0$-dimensional and $1$-dimensional orbits.
It turns out that the cosupports of $1$-dimensional
orbits are unions of two perfect matchings and
the cosupports of $0$-dimensional orbits are unions
of three perfect matchings. This allows us to reconstruct
the toric diagram of the toric $3$-Calabi-Yau variety $\lM_\te$. 

We will show that with any finite abelian group $G\sb\SL_3(\cC)$
we can associate a brane tiling. The corresponding quiver
potential algebra $\cC Q/(\dd W)$ is a non-commutative
crepant resolution of the quotient singularity $\cC^3/G$.
It is known that the Hilbert scheme $\Hilb^G(\cC^3)$
of $G$-clusters in $\cC^3$ is a crepant resolution
of $\cC^3/G$ \cite{Nakamura1}. This Hilbert scheme is
isomorphic to $\lM_\te$ for certain $\te$ (see Remark \ref{ghilb}).
It was shown by Nakamura \cite{Nakamura1} that
$\Hilb^G(\cC^3)$ is a toric variety. He also described
the corresponding fan.
According to the results mentioned above, we can
describe the toric diagram of $\lM_\te$ for any generic $\te$
by using the perfect matchings of the brane tiling.
A different algorithm to determine this toric diagram,
by computing the vertices of the polyhedron defining $\lM_\te$,
was proposed in \cite{Craw2}.

The paper is organized as follows: 
In Section \ref{prelim} we gather preliminary material
on brane tilings, quiver potential algebras,
consistency conditions, and Calabi-Yau property.
In Section \ref{sec:non-comm} we study some properties
of the quiver potential algebra induced by the brane
tiling and prove, in particular, that it is a non-commutative
resolution.
In Section \ref{sec:construction} we give a toric
description of the moduli spaces $\lM_\te$. We study
the orbits of $\lM_\te$ and give an explicit description
of its toric diagram.
In Section \ref{sec:orbifolds} we relate the McKay
correspondence for finite abelian groups $G\sb\SL_3(\cC)$
with brane tilings.

In the subsequent paper, joint with Martin Bender, we will
give a toric description of tilting bundles on
the crepant resolutions. This result gives a proof of
the conjecture
of Hanany, Herzog and Vegh \cite{HHV} and of a version
of the conjecture of Aspinwall~\cite{Aspinwall}.

I would like to thank Markus Reineke
for many useful discussions. 
I would also like to thank Igor Burban,
Alastair Craw, and Victor Ginzburg for many useful comments.

%% file: s01.tex
\section{Preliminaries}\label{prelim}

\subsection{Brane tilings}\label{pre:brane}
\begin{dfn}
A brane tiling is a bipartite graph $G=(G_0^\pm,G_1)$
together with an embedding of the corresponding CW-complex 
into the real two-dimensional torus $T$ so that the complement
$T\ms G$ consists of simply-connected components.
We call the elements of $G_0^+$ (resp.\ $G_0^-$)
white vertices (resp.\ black vertices).
We identify homotopy equivalent embeddings.
The set of connected components of $T\ms G$ is denoted by 
$G_2$ and is called the set of faces of $G$.
\end{dfn}

We define a quiver $Q=(Q_0,Q_1)$ dual to the graph $G$
as follows. The set of vertices $Q_0$ is $G_2$, the set
of arrows $Q_1$ is $G_1$. For any arrow $a\in Q_1$ we define
its endpoints to be the polygons in $G_2$
adjacent to $a$.
The direction of $a$ is chosen in such a way that
the white vertex is on the right of $a$.
For any arrow $a\in Q_1$, we define $s(a),t(a)\in Q_0$ to be
its source and target.
The CW-complex corresponding to $Q$ is automatically
embedded in $T$. 
The set of connected components of the complement,
called the set of faces of $Q$,
will be denoted by $Q_2$. It can be identified
with $G_0$. There is a decomposition $Q_2=Q_2^+\cup Q_2^-$
corresponding to the decomposition $G_0=G_0^+\cup G_0 ^-$.
It follows from our definition that the arrows
of the faces from $Q_2^+$ go clockwise and the arrows 
of the faces from $Q_2^-$ go anti-clockwise.

For any face $F\in Q_2$, we denote by $w_F$
the necklace (equivalence class of cycles in $Q$
up to shift) obtained by going along the arrows of $F$.
We define the potential of $Q$ (see e.g. \cite{Ginz1,Bock1}) by
$$W=\sum_{F\in Q_2^+}w_F-\sum_{F\in Q_2^-}w_F.$$
For any cycle $u=a_1\dots a_n$ in $Q$ and for any arrow $a\in Q_1$,
we define the differential
$$\frac{\dd u}{\dd a}=\sum_{i:a_i=a}a_{i+1}\dots a_na_1\dots a_{i-1}\in\cC Q.$$
Extending the differential by linearity, we get $\dd W/\dd a\in\cC Q$.

\begin{dfn}
Define a two-sided ideal $(\dd W)\sb\cC Q$
to be generated by $\dd W/\dd a$, $a\in Q_1$.
Define a quiver potential algebra
to be the algebra $\cC Q/(\dd W)$.
\end{dfn}

\subsection{Some groups related to brane tilings}\label{some groups}
Consider a complex of abelian groups
$$\cZ^{Q_2}\arr^{d_2}\cZ^{Q_1}\arr^{d_1}\cZ^{Q_0},$$
where $d_2(F)=\sum_{a\in F}a$, $F\in Q_2$ and $d_1(a)=t(a)-s(a)$
for any arrow $a\in Q$. 
Its homology groups coincide with the homology
groups of the $2$-dimensional torus containing $Q$.
Following \cite{MR2}, we define
$$\La=\cZ^{Q_1}/\angs{d_2(F)-d_2(G)}{F,G\in Q_2}.$$
Equivalently, \La is given by a cocartesian left upper square
\begin{diagram}
\cZ^{Q_2}&\rTo^{d_2}&\cZ^{Q_1}&\rTo^{d_1}&\cZ^{Q_0}\\
\dTo&&\dTo^{\wt}&\ruDotsto_d\\
\cZ&\rTo^{\om_\La}&\La\\
\dDotsto^{\om_M}&\ruTo_i\\
M
\end{diagram}
where the left arrow is given by $F\mto1,\ F\in Q_2$.
It is proved in \cite[Lemma 3.3]{MR2}, under the
condition that $G$ has at least one perfect matching,
that $\La$ is a free abelian group and the map
$\om_\La:\cZ\arr\La$ is injective.
There exists a unique map $d:\La\arr\cZ^{Q_0}$ making the right
triangle commutative. Note that $d\om_\La=0$. 
We put $M=\ker(d)$. Then there exists a
unique map $\om_M:\cZ\arr M$ making the lower triangle commutative.

Let $\b=\ker(\cZ^{Q_0}\arr\cZ)$, where the map is given by
$i\mto1$, $i\in Q_0$. This group is generated by the elements
of the form $j-i$, where $i,j\in Q_0$. This implies
that $\b=\im d_1=\im d$, as the quiver is connected.
We have then an exact sequence
$$0\arr M\arr^i\La\arr^d \b\arr 0.$$
Let us compute the ranks of the groups in this
exact sequence. It is clear that $\rk B=\#Q_0-1$.
We have
$$\rk\La=\rk(\coker d_2)+1=\rk(\ker d_1) +3=\rk\cZ^{Q_0}-\rk(\coker d_1)+3=\#Q_0+2.$$
This implies that $\rk M=3$.

\subsection{Weights and equivalence relations}\label{sec:weights}
We define a weak path in $Q$ to be a path consisting of
arrows of $Q$ and their inverses (for any arrow $a$ we
identify $aa\inv$ and $a\inv a$ with trivial paths).
For any weak path $u$, we define its content $|u|\in\cZ^{Q_1}$
by counting every arrow of $u$ with appropriate sign.
We define the weight of $u$ by $\wt(u)=\wt(|u|)\in\La$.

Let $\a=\cC Q/(\dd W)$ and let $\a'$ be obtained from $\a$ by inverting
all arrows.
We say that two paths in $Q$ are equivalent
if they are equal in $\a$. We say that
two weak paths are weakly equivalent if they are equal
in $\a'$ (cf.\ \cite[Section 4]{MR2}).
It is proved in \cite[Prop.\ 4.8]{MR2} 
under the condition that $G$ has at least one perfect matching

\begin{prp}\label{prp:weight equiv}
Two weak paths in $Q$ having the same start points
are weakly equivalent if and only if they have
the same weights.
\end{prp}

For any face $F\in Q_2$ and for any vertex $i\in F$ we define
$\om_{i,F}$ to be a cycle starting at $i$ and going along $F$.
It is proved in \cite{MR2} that $\om_{i,F}\sim\om_{i,G}$ if $i\in F\cap G$.
We denote the corresponding equivalence class by $\om_i$.
We denote its weight by $\ub\om$.

Let $\pi:\tl T\arr T$ be the universal cover of the torus and let
$\tl Q$ be the preimage of $Q$. Then $\tl Q$ is a periodic quiver.
We can define equivalence (resp.\ weak equivalence) relation
on the set of paths (resp.\ weak paths) of $\tl Q$ in the same 
way as above (see \cite{MR2}). For any weak paths $u$ in $\tl Q$
we can define a weak path $\pi(u)$ in $Q$.
We define then the weight $\wt(u)\in\La$ to be the weight of $\pi(u)$.
Similarly to Proposition \ref{prp:weight equiv} we can prove
that two weak paths in $\tl Q$ having the same start points
are weakly equivalent if and only if they have
the same weights.

\subsection{Consistency conditions}\label{consistency}
Let $G$ be a brane tiling and let $(Q,W)$ be
the corresponding quiver potential.

\begin{dfn}\label{cons cond}
A brane tiling $G$ is called consistent (resp.\ geometrically
consistent) if there exists a map
$R:Q_1\arr(0,1]$ (resp.\ $R:Q_1\arr(0,1)$), called an R-charge,
that satisfies
\begin{equation}
\sum_{a\in F}R_a=2,\qquad F\in Q_2,
\label{eq:charge1}
\end{equation}
\begin{equation}
\sum_{a\ni i}(1-R_a)=2,\qquad i\in Q_0.
\label{eq:charge2}
\end{equation}
\end{dfn}

\begin{rmr}
By the Birkhoff-von Neumann Theorem
(see e.g.\ \cite[Corollary 8.6a]{Schrijver1})
consistency condition implies 
that the bipartite graph $G$ is non-degenerate,
that is, every edge of $G$ is contained in some perfect matching.
\end{rmr}

\begin{thr}[{\cite[Theorem 8.15]{Bock2}}]
A brane tiling is consistent if and only
if any two weakly equivalent paths in $Q$
are equivalent.
\end{thr}

\begin{rmr}
It was proved earlier in \cite[Lemma 5.3.1]{HHV} that
geometric consistency implies
that any two paths in $Q$ having the same start points
and the same weight are equivalent. This implies
that weakly equivalent paths are equivalent.
\end{rmr}

\begin{rmr}
We do not discuss in this paper consistency conditions
on brane tilings involving zig-zag paths
(see \cite[Section 3.4.2]{Broomhead1}, \cite[Def.\ 5.2]{Ishii2},
\cite[Theorem 8.12]{Bock2}).
These consistency conditions are equivalent to Definition \ref{cons cond}
by \cite[Theorem 8.12]{Bock2} and \cite[Section 5]{Ishii2}.
\end{rmr}

\begin{dfn}\label{dfn:minmal}
We say that a path $u:i\arr j$ in $Q$ is minimal if it is
not equivalent to $v\om_i$ for any path $v:i\arr j$.
\end{dfn}

\begin{prp}[{see \cite[Lemma 7.3]{Davison1}}]
Assume that the brane tiling is consistent.
Then for any minimal path $u:i\arr j$ in $Q$ there
exists an arrow $a\in Q_1$ such that $s(a)=j$ and $au$
is still minimal.
\end{prp}

\begin{rmr}
This property together with a consistency condition was
used in \cite{MR2} (see also \cite{Broomhead1}, \cite{Davison1})
to show that the quiver potential algebra
is a $3$-Calabi-Yau algebra (see Section \ref{CY}).
\end{rmr}

Let \lA be the set of perfect matchings on $G$. 
Any perfect matching $I\in\lA$ can be considered as
a subset of $Q_1$, so we can define a linear map
$\hi_I:\cZ^{Q_1}\arr\cZ$
$$\hi_I(a)=\begin{cases}
	1,&a\in I,\\
	0,&a\not\in I.
\end{cases}$$
Note that $\hi_I(d_2(F))=1$ for any face $F\in Q_2$,
so we can factor $\hi_I:\cZ^{Q_1}\arr\cZ$ through $\La$ and get $\hi_I:\La\arr\cZ$.
Thus $\hi_I\in\La\dual$ and we can consider $\ub\hi_I:=i^*\hi_I\in M\dual$.
We define a cone $\si\sb M\dual_\cQ$
to be generated by $\ub\hi_I$, $I\in\lA$.
The following result is proved in \cite[Prop.\ 6.5]{Ishii2}

\begin{prp}\label{prp:extremal}
Let $I$ be some perfect matching in the consistent brane
tiling. Then the following conditions
are equivalent
\begin{enumerate}
	\item The ray in $M\dual_\cQ$ generated by $\ub\hi_I$ is an extremal
	ray of $\si$.
	\item For any $J\in\lA$, $J\ne I$ we have $\ub\hi_I\ne \ub\hi_J$.
	\item The quiver $Q\ms I=(Q_0,Q_1\ms I)$ is strongly connected,
	i.e.\ for any vertices $i,j\in Q_0$ there exists a path from $i$ to
	$j$ in $Q\ms I$.
\end{enumerate}
A perfect matching satisfying these conditions is called an extremal
(or corner, or external) perfect matching.
\end{prp}

In this paper we will work only with geometrically consistent brane
tilings because we will need the following important result proved by
Broomhead \cite[Prop.~6.2]{Broomhead1}

\begin{prp}\label{prp:broom1}
Assume that the brane tiling is geometrically consistent.
Then, for any vertices $i,j\in\tl Q$, there exists a
path $u:i\arr j$ such that $\hi_I(u)=0$ for some
extremal perfect matching.
\end{prp}

\begin{rmr}
It is conjectured that the analogous statement also holds
for consistent brane tilings. All the results of our
paper can be then proved in this generality.
\end{rmr}

\subsection{Calabi-Yau property}\label{CY}
In this section $A$ will be a 
(left and right) noetherian algebra, finitely
generated over a field $k$.
We define its enveloping algebra $A^e=A\ts_k A\op$. Then
$A$ is a module over $A^e$ in a natural way.

\begin{dfn}
We say that $A$ has finite Hochschild dimension
if $A$ has a finite projective resolution as an $A^e$-module.
We say that $A$ is homologically smooth if, moreover, this
resolution can be chosen to consist only of finitely generated $A^e$-modules.
\end{dfn}

\begin{rmr}\label{rmr:finite hochschild}
If $A$ has finite Hochschild dimension then $A$ and $A\op$ have finite
global dimension \cite[Ch.9, Prop.\ 7.6]{Cartan1}. 
In particular, $A$ has finite injective dimension
as a module over $A$ and over $A\op$ (we say that $A$ is Gorenstein
in this case).
\end{rmr}

\begin{dfn}\label{rem:global dim}
An algebra $A$ is called a Calabi-Yau algebra
of dimension $d$ if $A$ is homologically smooth and
$$\RHom_{A^e}(A,A\ts A)\iso A[-d]$$
in the category $D^b(A^e)$.
\end{dfn}

\begin{dfn}[{\cite[Def.\ 8.1]{Bergh5}, \cite[Def.\ 5.1]{Yek1}}]
An object $K\in D^b(A^e)$ is called a rigid dualizing complex
if
\begin{enumerate}
	\item $K$ has finite injective dimension over $A$ and $A\op$.
	\item The cohomologies of $K$ are finitely generated over $A$ and $A\op$.
	\item Canonical morphisms $A\arr\RHom_{A}(K,K)$ and 
	$A\arr\RHom_{A\op}(K,K)$ are isomorphisms in $D^b(A^e)$.
	\item (Rigidity) $\RHom_{A^e}(A,K\ts K)\iso K$ in $D^b(A^e)$.
\end{enumerate}
\end{dfn}

\begin{prp}[{\cite[Prop.\ 8.2]{Bergh5}}]
Any two rigid dualizing complexes in $D^b(A^e)$ are isomorphic.
\end{prp}

We will denote the dualizing complex of $A$ (if it exists)
by $K_A$. The following result gives an explicit description
of $K_A$ under certain conditions.

\begin{prp}[{\cite[Prop.\ 5.13]{Yek1}}]
Assume that $A$ is Gorenstein and has a rigid dualizing complex $K_A$.
Then
$K_A\iso\RHom_A(\RHom_{A^e}(A,A^e),A).$
\end{prp}

\begin{lmm}
Assume that $A$ is a $d$-CY algebra.
Then $A[d]$ is a rigid dualizing complex.
\end{lmm}
\begin{proof}
$A$ has finite injective dimension over $A$ and $A\op$ by Remark \ref{rem:global dim}.
From the Calabi-Yau property we get
$$\RHom_{A_e}(A,A[d]\ts A[d])=\RHom_{A^e}(A,A\ts A)[2d]\iso A[-d][2d]=A[d].$$
It follows that $A[d]$ satisfies all the conditions
on the rigid dualizing complex.
\end{proof}

\begin{prp}[{see \cite[Prop.\ 5.9]{Yek1}}]
If $A$ is finite over its center
and is finitely generated over $k$
then $A$ has a rigid dualizing complex.
More precisely, if $S\arr A$ is a finite central morphism, where
$S$ is a commutative smooth algebra of dimension $n$, then
$K_A:=\RHom_S(A,\Om^n_{S/k}[n])$ is a rigid dualizing complex over $A$.
\end{prp}

\begin{rmr}
If $S$ a commutative smooth algebra of dimension $n$ over $k$
then $S$ has a rigid dualizing complex $K_S=\Om^d_{S/k}[n]$.
If $R$ is a commutative algebra of finite type over $k$,
then we can always find a finite morphism $S\arr R$ with
$S$ smooth. Then $R$ has a rigid
dualizing complex $K_R=\RHom_S(R,K_S)$. If $R$
is a Cohen-Macaulay algebra of dimension $d$
then $K_R$ is concentrated in degree $-d$ and the
module $\om_R:=K_R[-d]$ is a canonical module of $R$
(see \cite[Def.\ 3.3.16]{HerzogBruns}).
Note that the canonical module of $R$
is defined only up to tensoring with
a projective $R$-module of rank $1$
(i.e.\ invertible sheaf on $\Spec R$).
The rigid dualizing complex $K_R$ of $R$ is,
in contrast, uniquely determined.
Thus $K_R[-d]$ gives a canonical choice
of a canonical module. 
\end{rmr}

\begin{thr}[{see \cite[Theorem 7.2.14]{Ginz1}}]\label{thr:CY}
Let $A$ have finite Hochschild dimension and let
$R\sb A$ be a central subalgebra, such that
$A$ is finitely generated as a module over $R$ and
$R$ is a Cohen-Macaulay domain, equidimensional
of dimension $d$.
Then the following conditions are equivalent
\begin{enumerate}
	\item $A$ is a $d$-CY algebra.
	\item $A$ is a maximal Cohen-Macaulay module over $R$ and
	$A\iso\Hom_R(A,K_R[-d])$, where $K_R$ is a rigid dualizing
	complex of $R$.
	\item For any $X\in D^b(\mod A)$, $Y\in D^-(\mod A)$,
	we have (functorially)
	$$\RHom_A(X,Y[d])\iso D\RHom_A(Y,X),$$
	where $D:D^+(\mod R)\arr D^-(\mod R)$ is defined by
	$Z\mto \RHom_R(Z,K_R)$.
\end{enumerate}
\end{thr}

\begin{rmr}
The algebra $R$ is finitely generated over $k$ under the
conditions of the theorem. 
This follows from \cite[Prop.~7.8]{AtiyahMacdonald}
if $A$ is commutative. For the non-commutative $A$
the proof goes through the same lines.
\end{rmr}

\begin{lmm}
Let $R$ be a commutative algebra of finite type
over an algebraically closed field $k$.
Assume that $R$ is Gorenstein and is equidimensional
of dimension $d$.
Let $\mod_\fl R$ be the category of finite
length $R$-modules. Then the following contravariant endofunctors
on $\mod_\fl R$ are isomorphic.
\begin{enumerate}
	\item $M\mto\Hom_R(M,E_R)$, where $E_R=\oplus_{m\in\Specm R}E(R/m)$.
	\item $M\mto\Hom_k(M,k)$.
	\item $M\mto\Ext^d_R(M,R)$.
	\item $M\mto\RHom_R(M,K_R)$.
\end{enumerate}
\end{lmm}
\begin{proof}
According to \cite[Prop.\ 1.1]{Ooishi}, there exists
a unique (up to equivalence) contravariant, exact functor
$D:\mod_\fl R\arr\mod_\fl R$ such that $D(R/m)\iso R/m$
for every $m\in\Specm R$. 
The first and the second functors obviously satisfy these conditions
(we use here the assumption that $k$ is algebraically closed).
It follows from our assumptions that $\om_R:=K_R[-d]$
is an invertible sheaf on $\Spec R$.
Recall that $R$ has an injective resolution \cite[\S 1]{Bass1}
$$0\arr R\arr\bigoplus_{\height p=0}E(R/p)\arr
\bigoplus_{\height p=1}E(R/p)\arr\dots\arr
\bigoplus_{\height p=d}E(R/p)=E_R\arr0.$$
This implies that $\Ext^d_R(M,R)=\Hom_R(M,E)$ for any
$M\in\mod_\fl R$ and so the first and the third
functors are equivalent. This also implies
$$\RHom(M,K_R)=\RHom(M,\om_R[d])=\RHom(M,A[d])\ts\om_R=\Hom(M,E_R)\ts\om_R.$$
Therefore the fourth functor also satisfies the above
conditions.
\end{proof}

\begin{crl}
Under the conditions of the theorem, assume that $R$
is Gorenstein and $k$ is algebraically closed.
Then for any finite dimensional
$A$-modules $X,Y$ we have (functorially)
$$\Ext^{d-i}(X,Y)\iso \Hom_k(\Ext^i(Y,X),k),\qquad 0\le i\le d.$$
\end{crl}

%% file: s02.tex
\section{Non-commutative crepant resolution}\label{sec:non-comm}
Let $\a=\cC Q/(\dd W)$ be the quiver potential algebra associated to
a geometrically consistent brane tiling.
In this section we will show that its center $R=Z(\a)$
is a normal Gorenstein domain and $\a$ is a non-commutative
crepant resolution of $R$ in the sense of Van den Bergh.
Related questions are studied in \cite{Bock2}. 

Let $\La^+\sb\La$ be a semigroup generated by the weights
of arrows.
Define a cone $P\sb\La_\cQ$ by
$$P=\left\{\sum a_i\la_i\mid
a_i\in\cQ_{\ge0},\ \la_i\in\La^+\text{ for all }i\right\}.$$
The saturation of $\La^+$ is given by $\ub{\La^+}=P\cap\La$.
Recall that in Section \ref{some groups} we have constructed an
exact sequence of free abelian groups
$$0\arr M\arr^i\La\arr^d\b\arr0. $$
For any $i,j\in Q_0$, we define 
$$\La_{ij}=d\inv(j-i)\sb\La,\qquad \La^+_{ij}=\La_{ij}\cap\La^+.$$
We define $M^+=M\cap\La^+$.
The following result is a consequence of 
the Birkhoff-von Neumann Theorem 
(see e.g.\ \cite[Corollary 8.6a]{Schrijver1})

\begin{prp}
The dual cone $P\dual\sb\La\dual_\cQ$
is generated by $\hi_I$, $I\in\lA$ (see Section \ref{consistency}).
\end{prp}

\begin{rmr}\label{rmr:generators of cone}
It is proved in \cite[Lemma 2.3.4]{Broomhead1} that
the semigroup $P\dual\cap\La\dual$ is generated by
$\hi_I$, $I\in\lA$. Moreover, every $\hi_I$
generates a $1$-dimensional face of $P\dual$.
\end{rmr}

\begin{crl}\label{crl:saturation}
The cone $P$ equals
$$P=\sets{\la\in\La_\cQ}{\hi_I(\la)\ge0\ \forall I\in\lA}.$$
The saturation of $\La^+$ equals
$$\ub{\La^+}=P\cap\La=\sets{\la\in\La}{\hi_I(\la)\ge0\ \forall I\in\lA}.$$
\end{crl}

\begin{lmm}\label{lmm:strict}
If $u$ is a weak path in $Q$ such that
$\hi_I(u)\ge0$ for any extremal perfect matching
then $u$ is weakly equivalent to a strict path.
\end{lmm}
\begin{proof}
We consider $u$ as a weak path from $i$ to $j$ in the periodic quiver $\tl Q$.
Let $v:i\arr j$ be a strict path such that $\hi_I(v)=0$
for some perfect matching $I$ (see Prop.\ \ref{prp:broom1}).
Then $u=v\om^k$ for some $k\in \cZ$ (see \cite[Lemma 4.6]{MR2}) and we have
$0\le\hi_I(u)=\hi_I(v)+k=k$. This implies that $v\om^k$ is a strict path.
\end{proof}

\begin{crl}\label{crl:small}
A path $u$ in $\tl Q$ is minimal (see Def.~\ref{dfn:minmal})
if and only if $\hi_I(u)=0$ for some extremal perfect matching.
\end{crl}
\begin{proof}
Assume that $\hi_I(u)=0$ for some (not necessarily extremal)
perfect matching.
If $u=v\om^k$ for some path $v$ and some $k\ge0$
then $0=\hi_I(u)=\hi_I(v)+k\ge k$, so $k=0$. This
implies that $u$ is minimal.

Assume that $u$ is minimal. Then $u\om\inv$ is not equivalent
to any strict path. It follows from Lemma \ref{lmm:strict}
that $\hi_I(u\om\inv)<0$ for some extremal perfect matching.
This implies $\hi_I(u)=0$.
\end{proof}

\begin{crl}
For any $i,j\in Q_0$, we have $\La^+_{ij}=\La_{ij}\cap P$.
The semigroup $M^+$ is saturated.
\end{crl}
\begin{proof}
We have from Corollary \ref{crl:saturation}
$$\La_{ij}\cap P=\sets{\la\in\La_{ij}}{\hi_I(\la)\ge0\ \forall I\in\lA}.$$
For any $\la$ as above, we can find a weak path $u:i\arr j$
such that $\wt(u)=\la$. Then $\hi_I(u)\ge0\ \forall I\in\lA$
and it follows from Lemma \ref{lmm:strict} that $u$ is
weakly equivalent to a strict path, so $\la=\wt(u)\in\La^+$.
This proves the first statement. The second statement
follows from the first for $i=j$.
\end{proof}

\begin{rmr}
The semigroup $\La^+$ is not saturated 
(i.e.\ $\La^+\ne\La\cap P$) in general.
So the above corollary can be quite confusing,
as it says that $\La_{ij}\cap\La^+=\La_{ij}\cap P$.
Note that $\La$ is not the union of 
$\La_{ij},\ i,j\in Q_0$.
\end{rmr}

\begin{rmr}
For any $i,j\in Q_0$ we can identify $e_j\a e_i$
with a vector space $\cC[\La^+_{ij}]$
(we denote its basis elements by $e^\la$, $\la\in\La^+_{ij}$).
This is the content of the algebraic consistency condition
\cite[Definition 4.4.2]{Broomhead1} proved in \cite{Broomhead1}.
For any $\la\in\La_{ij}$, we denote by
$u^\la_{ij}$ an (equivalence class of a)
path from $i$ to $j$ having weight \la.
The above identification is given by
$u^\la_{ij}\mto e^\la$, $\la\in\La^+_{ij}$.
\end{rmr}

\begin{rmr}\label{CY:M+}
Let $P_M=P\cap M_\cQ$. Then the dual cone $P_M\dual\sb M_\cQ\dual$
is generated by $\ub\hi_I$ for extremal perfect matchings $I$
(see Proposition \ref{prp:extremal}). The elements $\ub\hi_I$
are contained in the hyperplane
$$\sets{y\in M_\cQ\dual}{\om_M^*(y)=1},$$
where $\om_M:\cZ\arr M$ was defined in Section \ref{some groups}.
This implies that $\Spec\cC[M^+]=\Spec\cC[P_M\cap M]$
is a toric Calabi-Yau variety.
Its toric diagram is defined as an intersection of the cone
$P_M\dual$ with the above hyperplane.
\end{rmr}

\begin{lmm}[{see \cite[Lemma 4.3.1]{Broomhead1}}]
The center of the quiver potential algebra
$\a=\cC Q/(\dd W)$ is isomorphic to $\cC[M^+]$.
\end{lmm}

Let $R=\cC[M^+]$. The inclusion $R\arr \a$ from the above
lemma is given
by $e^\la\mto\sum_{i\in Q_0}u^\la_{ii}$.
Note that every $e_j\a e_i$ is automatically an $R$-module.

\begin{lmm}
For any $i,j,k\in Q_0$, there is a canonical $R$-module isomorphism
$$\Hom_R(e_j\a e_i,e_k\a e_i)=e_k\a e_j.$$
\end{lmm}
\begin{proof}
There is a natural embedding of $e_k\a e_j$ in
$\Hom_R(e_j\a e_i,e_k\a e_i)$. We just have to prove
that it is a bijection. Let $f\in\Hom_R(e_j\a e_i,e_k\a e_i)$.
We can assume that there exists some weak path
$u:j\arr k$ such that $f(v)=uv$ for any path $v:i\arr j$.
We have to show that $u$ is equivalent to a strict path.
Let $u=\om^ku'$ for the minimal path $u'$ and $k\in\cZ$.
By Corollary \ref{crl:small} there exists an external
perfect matching $I$ such that $\hi_I(u')=0$.
We know that the quiver $Q\ms I$ is strongly connected
(see Proposition~\ref{prp:extremal}), so there exists a
strict path $u'':i\arr j$ such that $\hi_I(u'')=0$.
Then the path $f(u'')=\om^ku'u''$ is strict.
Therefore $0\le\hi_I(\om^ku'u'')=k+\hi_I(u')+\hi_I(u'')=k$.
This implies that $u'\om^k$ is a strict path.
\end{proof}

\begin{crl}
For any $i\in Q_0$, we have $\Hom_R(\a e_i,\a e_i)=\a$.
\end{crl}

\begin{prp}\label{prp:noncom resolution}
The algebra $R$ is a normal Gorenstein domain and
the algebra $\a=\cC Q/(\dd W)$ is its non-commutative
crepant resolution (see \cite[Definition 4.1]{Bergh1}).
\end{prp}
\begin{proof}
The algebra $R=\cC[M^+]$ is normal, as $M^+$ is saturated.
Let
$$P_M=P\cap M_\cQ.$$
Then $M^+=P_M\cap M$.
According to \cite[p.~126]{Oda1}, the algebra $\cC[M^+]$
is Gorenstein if and only if
$$\inn(P_M)\cap M=m+(P_M\cap M)$$
for some $m\in M$.
Let $\lA^e\sb\lA$ be the set of extremal perfect matchings.
We have
$$\inn(P_M)\cap M=\sets{\la\in P\cap M}{\hi_I(\la)>0\ \forall I\in\lA^e}.$$
This implies that, for any $\la\in \inn(P_M)\cap M$,
we have $\hi_I(\la-\ub\om)=\hi_I(\la)-1\ge0$, $I\in\lA^e$,
so $\la-\ub\om\in P_M\cap M$. It follows that
$$\inn(P_M)\cap M=\ub\om+(P_M\cap M)$$
and the algebra $\cC[M^+]$ is Gorenstein.

Let $i\in Q_0$. We have seen that $\a =\End_R(\a e_i)$.
The $R$-module $\a e_i$ is reflexive. Indeed, we have
$$\Hom_R(e_j\a e_i,R)=\Hom_R(e_j\a e_i,e_i\a e_i)=e_i\a e_j.$$
Taking again the dual, we see that $e_j\a e_i$ is reflexive
and therefore also $\a e_i$ is reflexive.
According to \cite[Lemma 4.2]{Bergh1}, we just have to show
that $\a$ has finite global dimension and $\a$ is a CM-module
over $R$.
It is proved in \cite{MR2} that $\a$ is a $3$-Calabi-Yau algebra.
This together with Remark \ref{rmr:finite hochschild}
implies that $\a$ has finite global dimension.
To show that $\a$ is a CM-module over $R$ we will apply
Theorem \ref{thr:CY}. To do this we have to show
that $\a$ is a finitely generated module over $R$.
It is enough to show that
$e_j\a e_i$ is finitely generated over $R$ for every $i,j\in Q_0$.
Let us choose some path $v:j\arr i$. Then
the map $e_j\a e_i\arr e_i\a e_i\iso R$, $u\mto vu$ is injective.
It follows that $e_j\a e_i$ is finitely generated over $R$,
as $R=\cC[M^+]$ is noetherian.
\end{proof}

%% file: s03.tex
\section{Construction of crepant resolutions}\label{sec:construction}
Let $(Q,W)$ be the quiver potential associated with
a geometrically consistent brane tiling and let $\a=\cC Q/(\dd W)$ be the
corresponding quiver potential algebra.
We denote by $\mod\a$ the category of finitely generated left $A$-modules.

Recall that $\La^+\sb\La$ is a semigroup generated by the weights of the
arrows from $Q_1$. It was shown in \cite{MR2} that 
$\La^+\cap(-\La^+)=\set0$. 
The semigroup $\La^+$ is not saturated and 
$\Spec\cC[\La^+]$ is not normal in general.

We will consider the moduli spaces of representations of $A$
having dimension vector $\al=(1,\dots,1)\in\cZ^{Q_0}$.
The space of $\cC Q$-representations of dimension \al is given
by
$$R(\cC Q,\al)(\cC)=\cC^{Q_1}.$$
Let $R(A,\al)\sb R(\cC Q,\al)$ be the subvariety of those representations
that satisfy relations induced by the potential $W$.
The structure ring of $R(\cC Q,\al)$ is a polynomial algebra
$\cC[\cN^{Q_1}]$. We will denote its natural basis
by $(e^\la)_{\la\in\cN^{Q_1}}$.
For any arrow $a\in Q_1$, let $F^\pm_a\in Q_2$ be the faces
containing $a$. Let $u^\pm_a$ be a path in $Q$ such that $au^\pm_a$
is a cycle along $F^\pm_a$. Then
$$\frac{\dd W}{\dd a}=u^+_a-u^-_a.$$
This implies that the structure ring of $R(A,\al)$ is
given by
$$\cC[\cN^{Q_1}]/(e^{|u_a^+|}-e^{|u^-_a|}\mid a\in Q_1),$$
where, for any path $u$, the vector $|u|\in\cZ^{Q_1}$ denotes its content
(see Section \ref{sec:weights}).
The natural surjective map $\cC[\cN^{Q_1}]\arr\cC[\La^+]$
can be factored
$$\cC[\cN^{Q_1}]/(e^{|u_a^+|}-e^{|u^-_a|}\mid a\in Q_1)\arr \cC[\La^+].$$
This implies that there is a closed embedding
$$\Spec\cC[\La^+]\arr R(A,\al).$$

\begin{rmr}
This map need not be an isomorphism.
We thank Alastair Craw for this remark.
\end{rmr}

\begin{prp}
Variety $\Spec\cC[\La^+]$ is an irreducible component of $R(A,\al)$
\end{prp}
\begin{proof}
Recall that we have a map $\wt:\cZ^{Q_1}\arr\La$ (see Section \ref{some groups}). 
We define ideals $I,J\sb\cC[\cN^{Q_1}]$ by the rule
$$I=\left(e^{\la^+}-e^{\la^-}\mid \la^\pm\in\cN^{Q_1},\ \la^+-\la^-\in\ker(\wt)\right).$$
$$J=\left(e^{|u^+_a|}-e^{|u^-_a|}\mid a\in Q_1\right).$$
Note that the elements $|u^+_a|-|u^-_a|$, $a\in Q_1$, generate the group $\ker(\wt)$
(\La was defined as a factor group of $\cZ^{Q_1}$ by the subgroup
generated by the above elements, see Section \ref{some groups}).
Ideal $I$ defines variety $\Spec\cC[\La^+]$ and ideal $J$
defines variety $R(A,\al)$.
The proof now literally repeats the proof of \cite[Theorem 3.10]{Craw2}.
\end{proof}

In order to construct the moduli spaces of left $\a$-modules
of dimension \al, we have to identify isomorphic
representations from $R(\a,\al)$ with each other. This is achieved by
taking GIT quotients with respect to 
the natural action of the group
$$\GL_\al(\cC)=\prod_{i\in Q_0}\GL_{\al_i}(\cC)$$
on $R(\a,\al)$ \cite{King1}. 
One can see that two representations from $R(\a,\al)$
are isomorphic if and only if they are contained in the same orbit.
In our case
$$\GL_\al(\cC)=(\cC^*)^{Q_0}=\Hom_\cZ(\cZ^{Q_0},\cC^*).$$
The diagonal $\cC^*\sb\GL_\al(\cC)$ acts trivially on
$R(\a,\al)$, so our action factors through
$$T_\b=\Hom_\cZ(\b,\cC^*),$$
where $\b=\ker(\cZ^{Q_0}\arr\cZ)$ was defined in Section \ref{some groups}.
Given an element $\te\in \b$, we say that a representation $X\in R(A,\al)$
is stable (resp.\ semistable) if for any proper nonzero subrepresentation
$Y\sb X$, we have $\te\cdot\dim Y>0$ (resp.\ $\te\cdot\dim Y\ge0$).
According to \cite{King1}, the moduli space of \te-semistable 
left $\a$-modules of dimension \al is given by the GIT quotient
$$\lM_\te(\a,\al)=R(\a,\al)\GIT_\te T_\b.$$


\begin{dfn}
An element $\te\in\b$ is called \al-generic if,
for any $0<\be<\al$ we have $\te\cdot\be\ne0$.
\end{dfn}

If $\te\in \b$ is \al-generic, then all \te-semistable $\a$-modules
of dimension \al are stable.  

\begin{prp}
Assume that $\te\in B$ is \al-generic. Then we
have 
$$\lM_\te(A,\al)=\Spec\cC[\La^+]\GIT_\te T_B=\Spec\cC[P\cap\La]\GIT_\te T_B.$$
\end{prp}
\begin{proof}
We know from \cite[Prop.~5.1]{Ishii1} that $\lM_\te$ is
irreducible.
We can find some irreducible component $Z\sb R(A,\al)$
such that $Z\GIT_\te T_B$ equals $\lM_\te=\lM_\te(A,\al)$.
As all points of $T_\La\sb R(A,\al)$ are \te-stable
(they correspond to simple modules),
we deduce that $T_\La\sb Z$.
Variety $T_\La$ is dense in $\Spec\cC[\La^+]$, so 
$\Spec\cC[\La^+]\sb Z$. By the previous proposition
this inclusion is an isomorphism. This
implies $\Spec\cC[\La^+]\GIT_\te T_B=\lM_\te$.
To prove the second equality, we note
that $\Spec\cC[P\cap\La]\GIT_\te T_B$ is a normalization
of $\Spec\cC[\La^+]\GIT_\te T_B=\lM_\te$
and $\lM_\te$ is smooth by \cite[Prop.~5.1]{Ishii1}.
\end{proof}

If $\te\in B$ is \al-generic then
there exists a universal vector bundle $\lU$ on $\lM_\te=\lM_\te(\a,\al)$.

\begin{thr}\label{der equiv}
Let $\te\in\b$ be \al-generic.
Then $\lM_\te$ is smooth, the natural map $\lM_\te\arr \Spec Z(\a)$ is a crepant
resolution, and there is a pair of inverse equivalences of derived categories
$$\Phi:D^b(\coh \lM_\te)\arr D^b(\mod A\op),\quad F\mto R\Ga(F\ts_{\lM_\te}^L\lU^*),$$
$$\Psi:D^b(\mod \a\op)\arr D^b(\coh \lM_\te),\quad M\mto M\ts_A^L\lU.$$
\end{thr}
\begin{proof}
This follows from \cite[Theorem 6.3.1]{Bergh1} and the fact that $\a$ is
a non-commu\-ta\-ti\-ve crepant resolution of $Z(\a)$ (see Proposition
\ref{prp:noncom resolution}).
We should just note that Van den Bergh considers
the moduli spaces of right $A$-modules while we consider
the moduli spaces of left $A$-modules.
\end{proof}

\begin{rmr}
This result was proved directly by Ishii and Ueda \cite{Ishii1,Ishii2}.
\end{rmr}

According to \cite{Thaddeus1}, variety
$$\lM_\te(\a,\al)=\Spec\cC[P\cap\La]\GIT_\te T_\b$$
is a toric variety
endowed with an action of the torus $T_M=T_\La/T_\b$.
We are going to describe $T_M$-orbits of $\lM_\te(A,\al)$.
To do this we will describe the $T_\La$-orbits
of $R(\a,\al)$. It turns out that orbits corresponding
to indecomposable modules are parametrized by their supports.

\begin{dfn}
For any subset $I\sb Q_1$ we define representation
$x_I=(x_{I,a})_{a\in Q_1}\in R(\cC Q,\al)$ by the rule
$$x_{I,a}=
\begin{cases}
0,&a\in I\\
1,&a\not\in I
\end{cases}$$
We say that $I$ is $W$-compatible if $x_I$ is contained
in $R(\a,\al)$. We say that $I$ is \te-stable
(resp.\ semistable) for $\te\in \b$, if $x_I$ is \te-stable
(resp.\ semistable). We say that $I$ is indecomposable if
$x_I$ is indecomposable. 
\end{dfn}

\begin{rmr}
Note that any perfect matching $I\sb Q_1$ is $W$-compatible.
If $I$ is an extremal perfect matching then $Q\ms I$ is
a strongly connected quiver. This implies that $x_I$
is a simple representation and, in particular, \te-stable
for any $\te\in\b$.
\end{rmr}

\begin{dfn}
For any representation $x=(x_a)_{a\in Q_1}\in R(\cC Q,\al)$
we define its cosupport
$$I_x=\sets{a\in Q_1}{x_a=0}.$$
All representations in the same $(\cC^*)^{Q_1}$-orbit
of $R(\cC Q,\al)$ have the same cosupport.
\end{dfn}

\begin{lmm}
Let $x\in R(\a,\al)$. Then $I_x$ is $W$-compatible.
Representation $x$ is \te-stable (resp.\ \te-semistable,
resp.\ indecomposable) if and only if $I_x$ is
\te-stable (resp. \te-semistable, resp.\ indecomposable).
\end{lmm}

\begin{dfn}
Let $I\sb Q_1$.
We consider it as a subgraph of the bipartite graph
$G=(G_0^\pm,G_1)$ dual to $Q$ in the 
two-dimensional torus $T$.
A connected component of $I$ is called a big component
of $I$ if it contains more then one edge.
Otherwise it is called a small component of $I$.
\end{dfn}

\begin{prp}
Two indecomposable representations $x,y\in R(\a,\al)$
are contained in the same $T_\La$-orbit if and only
if they have the same support.
\end{prp}
\begin{proof}
It is clear that any two representations in the same 
$T_\La$-orbit have equal supports. Let us prove the converse.
So let $x,y\in R(\a,\al)$ be such that $I_x=I_y=:I$.
We define a new Q-representation $z$ by the rule
$$z_a=\begin{cases}
	0,&a\in I\\
	x_a\inv y_a,&a\not\in I
\end{cases}$$
It is clear that $z$ is also an $\a$-representation.
It is also indecomposable, as this is a property
of the support. 
We claim that we can extend $(z_a)_{a\in Q_1\ms I}$ to the element
in $(\cC^*)^{Q_1}$ that is still an $\a$-representation.
Such an element will be automatically contained in $T_\La$
and will map $x$ to $y$, so that both elements will be in the same
$T_\La$-orbit.

We may suppose that $I\ne\emptyset$, as otherwise our
claim is automatically satisfied.
By the $W$-relations the elements $\prod_{a\in F}z_a$
are independent of $F\in Q_2$ and therefore are all zero.
It follows that every face intersects $I$ non-trivially.
If every face intersects $I$ in precisely one element,
then $I$ is a perfect matching. In this case $z$ can be
obviously extended to $(\cC^*)^{Q_1}$ and we are done.

So let us assume that $I$ is not a perfect matching.
We consider $I$ as a subgraph of the bipartite graph
$G=(G_0^\pm,G_1)$ dual to $Q$ in the 
two-dimensional torus $T$.
Graph $I$ can have many connected
components.
Our assumptions imply that the set of vertices of $I$ equals
the set of vertices of $G$,
that every vertex (i.e.\ face in $Q$) has valence at least one,
and that there is at least one vertex 
having valence $\ge2$.
It follows from $W$-relations that any vertex connected
to a vertex of valence $\ge2$ also has valence $\ge2$.
So, if a connected component of $I$ is big then
all its vertices have valence $\ge2$.
We claim that there is precisely one big component.

First, let us note that the complement of $I$ in the torus
is connected as $z$ is indecomposable. 
One can easily see that the complement of two
non-intersecting loops in the torus always has
at least two connected components. Every
big component of $I$ has loops, so the
existence of two big components
would imply that the complement of $I$ is not
connected. This proves our claim
that there exists just one big component in $I$.

We will extend $(z_a)_{a\in Q_1\ms I}$ to the element
in $(\cC^*)^{Q_1}$ in such a way that products along
the faces (elements in $Q_2$) are all equal one. The choice
for the arrows of small components is clear. Let $J\sb I$
denotes the big component. For any $F\in G_0=Q_2$
we define 
$$z_F=\prod_{a\in F\ms I}z_a^{\eps(F)},$$
where $\eps(F)=\pm1$ for $F\in G_0^\pm$. 
It follows from the $W$-relations that
$$\prod_{F\in J_0}z_F=1,$$
where $J_0\sb G_0$ (resp.\ $J_1\sb G_1$)denotes the set of vertices of $J$
(resp.\ the set of edges of $J$).
More generally, we consider an arbitrary abelian group
\Ga and a sequence
$$\Ga^{J_1}\arr^d\Ga^{J_0}\arr^p \Ga,$$
where $d(\ga a)=\ga F^+_a-\ga F^-_a$ with $F_a^\pm\in J_0^\pm$ 
incident with $a\in J_1$, and $p(\ga F)=\ga$ for $F\in J_0$, $\ga\in\Ga$.
The kernel of $p$ is generated by the elements of the form
$\ga F^+_a-\ga F^-_a$, $a\in J_1,\ga\in\Ga$, as $J$ is connected.
So the above sequence is exact.
In our situation this means that there exists
$(t_a)_{a\in J_1}\in(\cC^*)^{J_1}$ such that,
for any $F\in J_0\sb Q_2$, we have
$$\prod_{a\in F\cap J_1}t_a^{\eps(F)}=z_F.$$
The elements $t_a\inv,\ a\in J_1$ give then the desired
extension.
\end{proof}


\begin{rmr}
We have shown that the cosupport of an indecomposable $A$-module
of dimension \al can have at most one big component.
Every vertex of a big component has valence at least $2$.
\end{rmr}

\begin{crl}
There is a bijection between the $T_\La$-orbits of
indecomposable representations in $R(\a,\al)$
and $W$-compatible, indecomposable subsets in $Q_1$.
\end{crl}

Let now $\te\in \b$ be \al-generic. Then all 
representations from $\lM_\te=\lM_\te(\a,\al)$ are
stable and in particular indecomposable.
It follows that $T_\b$-orbits of $\lM_\te$
are parametrized by $W$-compatible, \te-stable
subsets in $Q_1$. We are going to give a precise
description of these subsets according to the dimension
of the corresponding orbits. For any $W$-compatible,
\te-stable set $I\sb Q_1$, we denote by $\si_I$ the cone
of the fan of $\lM_\te$ corresponding to the orbit defined by $I$.

\begin{prp}
Let $x\in\lM_\te$ and let $O_x\sb\lM_\te$ be its $T_M$-orbit.
Then
\begin{enumerate}
	\item $\dim O_x=3$ if and only if $I_x=\emptyset$.
	\item $\dim O_x=2$ if and only if $I_x$ is a perfect matching.
	\item $\dim O_x=1$ if and only if the big component of $I_x$ is a cycle.
	\item $\dim O_x=0$ if and only if the big component of $I_x$ has two
	trivalent vertices of different colors and all other vertices
	of valence $2$.
\end{enumerate}
\end{prp}
\begin{proof}
The first statement is trivial.
The last statement is just a translation of \cite[Lemma 4.5]{Ishii1}
to our language in the case of a consistent brane tiling.
It is proved there also that if \si is a cone
corresponding to the fixed point $x$ then $U_\si=\Spec\cC[\si\dual\cap M]$
is isomorphic to $\cC^3$.
This isomorphism can be described in the following way.
Let $v^\pm$ be the white and black vertices of valence $3$ in
the big component of $I_x$. Let $a^\pm_i,i=1,2,3$ be the
edges in the big component incident with $v^\pm$.
Then for any choice of values $t_i\in\cC,i=1,2,3$ for the edges
$a^+_i,i=1,2,3$ there exists the unique choice of
the values for the rest of the edges in $I_x$ such that the
corresponding quiver representation satisfies $W$-relations
(representation $x$ corresponds to zero values on the
edges $a^+_i,i=1,2,3$). This gives the required isomorphism
$U_\si\iso\cC^3$.

For any one-dimensional orbit $O_y$ there exists a fixed point
$x$ in its closure. Let \si be the cone corresponding
to the point $x$. Then $O_y\sb U_\si$ and using the
above identification $U_\si\iso\cC^3$, we
can write without loss of generality
$$O_y=\sets{(t,0,0)}{t\in\cC^*}.$$
Let
$$(a_{i,1},\dots,a_{i,2k_i+1}),\quad k_i\ge0,\quad i=1,2,3$$
be the chains of edges in the big component of $I_x$
that connect $v^+$ and $v^-$. We assume that $a_{i,1}=a^+_i$,
$i=1,2,3$.
Then representation $y$ satisfies
\begin{align*}
y_{a_{1,2k+1}}\ne0,	&\qquad 0\le k\le k_1\\
y_{a_{1,2k}}=0,		&\qquad 1\le k\le k_1\\
y_{a_{i,k}}=0,		&\qquad 1\le k\le 2k_i+1,\ i=2,3.
\end{align*}
This implies that the big component of $I_y$ consists
of the edges $a_{i,k}$, $1\le k\le 2k_i+1,i=2,3$. This
is a cycle.

For any two-dimensional orbit $O_y$ we again consider a fixed point
$x$ in its closure and the cone $\si$ corresponding to $x$.
We can write without loss of generality
$$O_y=\sets{(t_1,t_2,0)}{t_1,t_2\in\cC^*}.$$
We use the same notation as above for the chains connecting $v_+$
and $v_-$. Then representation $y$ satisfies
\begin{align*}
y_{a_{i,2k+1}}\ne0,	&\qquad 0\le k\le k_i,i=1,2\\
y_{a_{i,2k}}=0,		&\qquad 1\le k\le k_i,i=1,2\\
y_{a_{3,2k+1}}=0,	&\qquad 0\le k\le k_3\\
y_{a_{3,2k}}\ne0,		&\qquad 1\le k\le k_3.
\end{align*}
This implies that $I_y$ is a perfect matching.
\end{proof}

\begin{rmr}
The one-dimensional cones of the fan of $\lM_\te$ are generated
by $\ub\hi_I\in M\dual$, where $I\in\lA$ are \te-stable.
All these vectors are contained in the hyperplane
$$\sets{y\in M_\cQ\dual}{\om_M^*(y)=1},$$
where $\om_M:\cZ\arr M$ was defined in Section \ref{some groups}.
This implies that $\lM_\te$ is a toric Calabi-Yau variety
(cf.\ Remark \ref{CY:M+}).
Its toric diagram is defined as an intersection of the cones of 
the fan of $\lM_\te$
with the above hyperplane. This toric diagram is a triangulation
of the toric diagram of $\cC[M^+]$ (see Remark \ref{CY:M+}).
\end{rmr}

\begin{rmr}
For any non-extremal point of the toric diagram there exists
more than one perfect matching that is mapped to it 
(see Proposition \ref{prp:extremal}).
However, only one of these perfect matchings is \te-stable.
\end{rmr}

\begin{crl}
Let $x\in\lM_\te$ be a fixed point and let
$$(a_{i,1},\dots,a_{i,2k_i+1}),\qquad i=1,2,3$$
be the chains of edges connecting the trivalent points
of the big component of $I_x$. Then there are precisely
three $2$-dimensional orbits containing $x$ in their closures.
The corresponding perfect matchings are given by (for $i=1,2,3$)
$$I_x\ms\sets{a_{j,2k+1}}{0\le k\le k_j,j\ne i}\ms\sets{a_{i,2k}}{1\le k\le k_i}.$$
These are the only perfect matchings contained in $I_x$.
$I_x$ is a union of these perfect matchings.
\end{crl}

\begin{crl}
Let $O_x\sb\lM_\te$ be a one-dimensional orbit and let
$$(a_1,\dots,a_{2k})$$
be the big component of $I_x$ which is a cycle. 
Then there are precisely two $2$-dimensional orbits containing
$O_x$ in their closures.
The corresponding perfect matchings are given by
$$I_x\ms\set{a_1,a_3,\dots,a_{2k-1}},\qquad I_x\ms\set{a_2,a_4,\dots,a_{2k}}.$$
These are the only perfect matchings contained in $I_x$.
$I_x$ is a union of these perfect matchings.
\end{crl}

We have now a simple algorithm to construct a fan of $\lM_\te$ for
\al-generic $\te\in \b$. 
Construction of such fan is equivalent to the
triangulation of the toric diagram.
We find first those perfect matchings that are \te-stable.
Then we test for all the pairs of these perfect
matchings if their union is \te-stable
(note that the union of any two perfect matchings is
automatically $W$-compatible). 
In this way we get all $1$-dimensional and $2$-dimensional
cones of the required fan.
The $3$-dimensional cones are constructed then automatically.
The triangulation of the toric diagram
is uniquely determined from this data.

%% file: s04.tex
\section{Orbifolds and brane tilings}\label{sec:orbifolds}
In this section, for any finite abelian subgroup $G\sb\SL_3(\cC)$,
we construct certain brane tiling.
The underlying quiver will coincide with the McKay quiver
of the $G$-representation $\cC^3$.
It turns out that the corresponding quiver potential algebra
is Morita-equivalent (and even isomorphic) 
to the skew product $\cC[x,y,z]\rtimes G$ 
and it is therefore the most natural candidate for the non-commutative
crepant resolution of the orbifold singularity $\cC^3/G$.

\subsection{McKay quiver}
Let $G$ be an arbitrary finite group and let $V$ be its finite-dimensional
representation over \cC. We denote by $\hat G$ the set of all irreducible
$G$-representations. 

\begin{rmr}
If $G$ is abelian, then $\hat G$ has a group
structure induced by tensor products and $\hat G$ can be identified
with $\Hom_\cZ(G,\cC^*)$. It is called the group
of characters of $G$.
This group is non-canonically isomorphic to $G$.
\end{rmr}

Recall the definition of the McKay quiver $Q$ of a $G$-representation $V$.
Its set of vertices is given by $Q_0=\hat G$. 
The set of arrows from $\si\in\hat G$ to $\rho\in\hat G$ is given
by some fixed basis of
$$\Hom_G(\si,\rho\ts V).$$
This quiver can be used to describe left modules
over a skew-group algebra $S(V\dual)\rtimes G$.
To define such module structure on a vector
space $M$ it is enough to endow $M$ with a structure of
a $G$-representation and to give a $G$-equivariant linear map $V\dual\ts M\arr M$
(corresponding to the multiplication) that satisfies certain axioms
(see e.g.\ \cite[Section 3]{Ito-Nakajima1}).
Let us show that this data defines a representation of the McKay quiver.
We decompose $M=\oplus_{\rho\in\hat G}M_\rho\ts\rho$,
where $M_\rho$
are some vector spaces. We put this vector spaces at the vertices
of the McKay quiver (recall that $Q_0=\hat G$). We have
\begin{multline*}
\Hom_G(V\dual\ts M,M)
=\oplus_{\si,\rho}\Hom_G(V\dual\ts\si,\rho)\ts\Hom(M_\si,M_\rho)\\
=\oplus_{\si,\rho}\Hom_G(\si,\rho\ts V)\ts\Hom(M_\si,M_\rho).
\end{multline*}
This means that for any arrow from $\si\in\hat G$ to $\rho\in\hat G$,
we have a linear map $M_\si\arr M_\rho$. This gives us a required quiver
representation. In this way we get a full and faithful functor
$$\mod (S(V\dual)\rtimes G)\arr \mod\cC Q.$$
To get an equivalence of
categories, we have to impose certain relations in
the path algebra $\cC Q$ \cite[Section 3]{Ito-Nakajima1}.

Let now $G$ be a finite abelian subgroup $G\sb\SL_3(\cC)$.
For any character $\rho\in\hat G$ we denote the corresponding
one-dimensional $G$-representation also by \rho.
We can decompose the $G$-representation $V=\cC^3$ (induced by
the inclusion $G\sb\SL(\cC^3)$)
$$V=\rho_1\oplus\rho_2\oplus\rho_3,$$
where $\rho_1,\rho_2,\rho_3\in\hat G$. Note that $\rho_1\rho_2\rho_3=1$,
as $G\sb\SL_3(\cC)$.
Then
$$\rho\ts V\iso\rho\rho_1\oplus\rho\rho_2\oplus\rho\rho_3$$
and $\Hom_G(\si,\rho\ts V)$ can be nonzero only if $\si=\rho\rho_i$
for some $i=1,2,3$. It follows that any vertex $\rho\in Q_0=\hat G$
has three ingoing arrows
$$a_i^\rho:\rho\rho_i\arr\rho,\qquad i=1,2,3.$$
Sometimes we will omit the upper index if the ingoing or outgoing vertex 
is known.
Now we describe the set of faces $Q_2$ corresponding to some brane tiling.
All faces will contain just three arrows. For any
vertex $\rho\in\hat G$ and any permutation $\pi\in S_3$, 
we consider the face
$$\rho\arr^{a_{\pi(3)}}
\rho\rho_{\pi(1)}\rho_{\pi(2)}
\arr^{a_{\pi(2)}}
\rho\rho_{\pi(1)}
\arr^{a_{\pi(1)}}\rho.$$
It is clear that every arrow is contained
in precisely two faces. This implies that if we glue
the faces along the common arrows we obtain some oriented compact
surface. The number of arrows equals $3\cdot\# G$ and the
number of faces equals $2\cdot\# G$. This implies that
the Euler number of our surface is zero and therefore
the surface is homotopic to a torus. So we obtain a
brane tiling. 
Let $W$ be the corresponding potential.

\begin{rmr}
We have used the fact that $Q$ is a connected quiver.
This follows from the fact that $\rho_1,\rho_2,\rho_3$
generate the whole group $\hat G$. Indeed, assume that
they generate some proper subgroup $\hat H\sb\hat G$.
Then the map $G\emb\SL_3(\cC)$ can be factored through
$G\arr H\arr\SL_3(\cC)$, where $H=\Hom_\cZ(\hat H,\cC^*)$.
This would imply that $G\arr \SL_3(\cC)$ is not injective.
\end{rmr}

The next result follows from 
\cite[Section 5.2]{sardoinfirri-1996} or
\cite[Section 3]{Ito-Nakajima1} or
\cite[Section 2]{Craw2}

\begin{prp}
There is an equivalence of categories
$$\mod (\cC[x,y,z]\rtimes G)\arr\mod \cC Q/(\dd W)$$.
\end{prp}

\begin{rmr}
It is proved in \cite[Proposition 2.8]{Craw3}
that the above algebras are actually 
isomorphic.
\end{rmr}

\begin{rmr}
It was shown in \cite[Section 4.4]{Ginz1} that for any
finite subgroup $G\sb\SL_3(\cC)$ one can endow the 
corresponding McKay quiver with a potential (depending
on some choices) in
such a way that the quiver potential algebra is Morita
equivalent to the skew group algebra. The coefficients
of the cycles in that potential are not always $\pm1$,
so it can not correspond to some brane tiling. However
one can show that if $G$ is abelian then one can make
such choices that the coefficients of the cycles
of the potential are $\pm1$ and this potential
is induced by the brane tiling constructed above.
\end{rmr}

\begin{rmr}
The periodic quiver of the above brane tiling coincides
with the periodic quiver corresponding to the brane
tiling of $\cC^3$. This implies that the constructed brane
tiling is always geometrically consistent.
\end{rmr}

\subsection{Toric realization of orbifolds}
Let $G\sb\SL_3(\cC)$ be a finite abelian group.
We have associated a quiver potential $(Q,W)$ and a triple
of characters $\rho_1,\rho_2,\rho_3\in\hat G$ with
such a group.
These characters generate $\hat G$.
So, we get a surjective map
$$p:M_0=\cZ^3\arr\hat G.$$
There is a map $\pi:\cZ^{Q_1}\arr M_0$, that maps an arrow
$a^\rho_i$ to the $i$-th canonical basis vector of $M_0=\cZ^3$.
It can be factored through $\pi:\La\arr M_0$.
Recall that we have defined $\b=\ker(\cZ^{Q_0}\arr\cZ)$.
We define a map $\nu:\b\arr\hat G$ to be the composition
$\b\emb\cZ^{Q_0}=\cZ^{\hat G}\arr\hat G$.
It follows from \cite[Lemma 10.5]{sardoinfirri-1996}

\begin{lmm}
The following diagram is cartesian and cocartesian
\begin{diagram}
\La&\rTo^d&\b\\
\dTo^\pi&&\dTo_\nu\\
M_0&\rTo^p&\hat G
\end{diagram}
\end{lmm}

\begin{crl}
We have a commutative diagram with short exact sequences in the
rows
\begin{diagram}
0&\rTo&M&\rTo&\La&\rTo^d&\b&\rTo&0\\
&&\dEq&&\dTo^\pi&&\dTo_\nu\\
0&\rTo&M&\rTo&M_0&\rTo^p&\hat G&\rTo&0
\end{diagram}
\end{crl}

This corollary allows us to interpret all the results
of \cite{sardoinfirri-1996}, \cite{Nakamura1} or \cite{Craw2}
in the context of brane tilings. 
For example, let $P_0\sb (M_0)_\cQ=\cQ^3$ be a cone generated
by the basis vectors. Then $\cC^3=\Spec\cC[P_0\cap M_0]$
and $\cC^3/G=\Spec\cC[P_0\cap M]$ by the general theory
of toric quotients \cite[Prop.\ 3.1]{Thaddeus1}
(see also \cite[Section 1]{Nakamura1}).
But the last scheme coincides with $\Spec Z(A)=\Spec \cC[P\cap M]$,
where $A=\cC Q/(\dd W)$ and $P\sb\La_\cQ$ is a cone defined in
Section \ref{sec:non-comm}.

\begin{rmr}\label{ghilb}
It is proved in \cite{Nakamura1} that $\cC^3/G$ has a
crepant resolution $\Hilb^G(\cC^3)$ -- the Hilbert scheme
of $G$-clusters in $\cC^3$ (see \cite{Nakamura1,Craw2}).
It follows from \cite[Prop.~5.2]{Craw2} that $\Hilb^G(\cC^3)$
can be identified with $\lM_\te(A,\al)$, where $A=\cC Q/(\dd W)$,
$\al=(1,\dots,1)\in\cZ^{Q_0}$ and $\te\in\b$ is such
that $\te_{\rho_0}<0$ and $\te_\rho>0$ for $\rho\ne\rho_0$
($\rho_0\in \hat G$ is a trivial representation).
\end{rmr}

\begin{rmr}
It was shown in \cite{Nakamura1} that $\Hilb^G(\cC^3)$
is a toric variety and it was proposed there a way
to compute the corresponding fan. 
In \cite{sardoinfirri-1996,Craw2} it was shown that 
$\lM_\te(A,\al)$ is a toric variety for \al-generic 
$\te\in \b$ and in \cite{Craw2} an algorithm was given 
to compute the corresponding fan. This 
algorithm consists in computing the vertices of the
polyhedron $P^\te\cap M_\cQ$ (this polyhedron defines a
toric variety $\lM_\te$, see \cite{Thaddeus1} for the
general results on toric varieties defined by polyhedra),
where $P^\te=P-\la$ for some $\la\in\La$ with $d(\la)=\te$. 
The results of the previous sections give us a new
way to compute this fan by using \te-stable
perfect matchings on the brane tiling constructed above.
\end{rmr}

%% file: example.tex
\subsection{Example}
In this example we will describe the toric diagram
of $\Hilb^G(\cC^3)$, where $G=\cZ_6$ and the
action on $\cC^3$ is given by $\frac16(1,2,3)$.
We have chosen this example as it was also
considered by Nakamura \cite{Nakamura1} using
completely different methods.

Let us make first some general remarks.
We have an exact sequence
$$0\arr M\arr M_0\arr \hat G\arr0.$$
Applying the functor $\Hom_\cZ(-,\cZ)$ we get an exact
sequence
$$0\arr\Hom(M_0,\cZ)\arr\Hom(M,\cZ)\arr\Ext^1(\hat G,\cZ)\arr0.$$
We claim that $\Ext^1(\hat G,\cZ)=G$.
The module \cZ over the ring \cZ has an injective resolution
$$0\arr\cZ\arr\cQ\arr\cQ/\cZ\arr0.$$
Using it we get 
$$\Ext^1(\hat G,\cZ)=\Hom(\hat G,\cQ/\cZ)=\Hom(\hat G,\cC^*)=G,$$
where we have used an inclusion $\cQ/\cZ\arr\cC^*$, $x\mto e^{2\pi ix}$.
This implies that there is an exact sequence
$$0\arr M_0\dual\arr M\dual\arr G\arr0.$$
Let $n=\# G$.
Then, for any $f\in M\dual$, we have 
$nf\in M_0\dual$. In particular,
for any perfect matching $I\in\lA$, we have
$n\ub\hi_I\in M_0\dual$.
To determine this function, we have to find its
values on the basis elements $e_i\in M_0$, $i=1,2,3$.
Given an arrow $a^\rho_i$, $\rho\in\hat G$, $i=1,2,3$,
we say that it has type $i$. 
We have $n\ub\hi_I(e_i)=\ub\hi_I(n e_i)$,
so to find this value, we have to evaluate
$\hi_I$ on any path consisting of $n$ arrows of type $i$.
Such a path will not necessarily contain only
pairwise different arrows
(there are precisely $n$ different arrows
of type $i$ in quiver $Q$). But one can
easily show that the value of $\hi_I$ on such
a path equals the number of arrows of type $i$
in the perfect matching $I$.
This gives us an easy way to determine the
elements $\ub\hi_I\in (M_0)_\cQ\dual=M_\cQ\dual$.

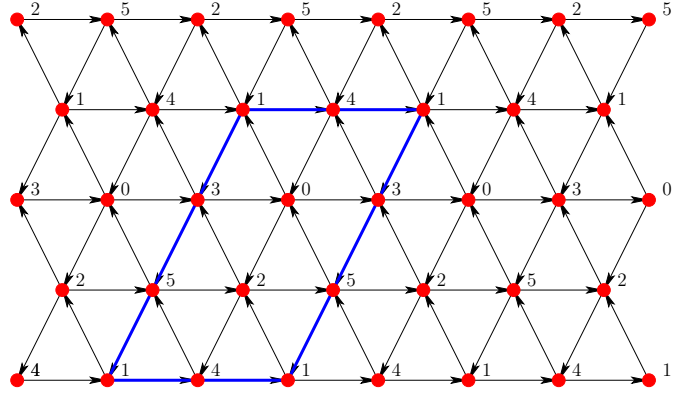
\begin{figure}[h]
\scalebox{0.6}{\input{Z6-123.pstex_t}}
\caption{Periodic quiver with a fundamental domain}
\label{fig:periodic}
\end{figure}

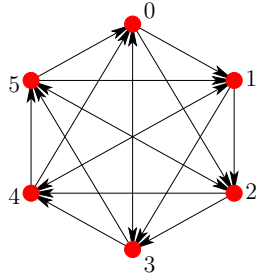
\begin{figure}[h]
\scalebox{0.75}{\input{Z6-123-mckay.pstex_t}}
\caption{McKay quiver}
\label{fig:mckay}
\end{figure}

Let us return now to the case $G=\cZ^6$ with an action $\frac16(1,2,3)$.
We depict first the periodic quiver and the fundamental
domain corresponding to the brane tiling constructed earlier,
see Figure \ref{fig:periodic}.
The corresponding quiver on the torus 
is a McKay quiver of 
the $G$-representation $\cC^3$, see Figure \ref{fig:mckay}

We will denote the arrow from vertex $i$ to
vertex $j$ by $ij$. The type of such arrow
equals $j-i$ if $j>i$ and $j-i+6$ otherwise.
The list of all perfect matchings of the brane tiling
is given in the Table \ref{table:pm}. Every perfect
matching is written as a set of arrows from $Q$.
\begin{table}[hp]
\begin{tabular}{|c|c|c|}
\hline
$N$ &$I$& $6\ub\hi_I$\\
\hline
1&$34, 01, 02, 35, 24, 51 $&$ (2, 4, 0)$\\
2&$34, 01, 02, 35, 45, 12 $&$ (4, 2, 0)$\\
3&$34, 01, 23, 50, 24, 51 $&$ (4, 2, 0)$\\
4&$34, 01, 23, 50, 45, 12 $&$ (6, 0, 0)$\\ 
5&$34, 14, 30, 50, 52, 12 $&$ (3, 0, 3)$\\ 
6&$34, 14, 30, 35, 24, 25 $&$ (1, 2, 3)$\\
7&$13, 40, 02, 35, 24, 51 $&$ (0, 6, 0)$\\ 
8&$13, 40, 02, 35, 45, 12 $&$ (2, 4, 0)$\\
9&$13, 40, 23, 50, 24, 51 $&$ (2, 4, 0)$\\ 
10&$13, 40, 23, 50, 45, 12 $&$ (4, 2, 0)$\\ 
11&$13, 14, 02, 03, 52, 12 $&$ (1, 2, 3)$\\
12&$13, 14, 23, 03, 24, 25 $&$ (1, 2, 3)$\\
13&$41, 40, 30, 50, 52, 51 $&$ (1, 2, 3)$\\ 
14&$41, 40, 30, 35, 45, 25 $&$ (1, 2, 3)$\\
15&$41, 03, 01, 02, 52, 51 $&$ (1, 2, 3)$\\
16&$41, 03, 01, 23, 45, 25 $&$ (3, 0, 3)$\\
17&$41, 03, 14, 52, 25, 30 $&$ (0, 0, 6)$\\ 
\hline
\end{tabular}
\bigskip
\caption{Perfect matchings and their coordinates in $M_\cQ\dual$}
\label{table:pm}
\end{table}
\begin{rmr}
We see that the elements $\frac16(1,2,3)$, $\frac16(2,4,0)$, etc.
are contained in $M\dual$.
Let $g\in G=\cZ_6$ be the image of $\frac16(1,2,3)$ with respect to
the canonical map $M\dual\arr G$. This is a generator of $G$.
We can see that $\frac16(2,4,0)\mto g^2$,
$\frac16(3,0,3)\mto g^3$, $\frac16(4,2,0)\mto g^4$ (cf.\ \cite{Nakamura1}).
These elements of $M\dual$ will be sometimes denoted by their images in $G$.
The perfect matchings $I$ such that $\ub\hi_I=e_i$, $i=1,2,3$,
are precisely the extremal perfect matchings. 
\end{rmr}

We consider now a stability $\te\in\b$ that corresponds to
the resolution $\Hilb^G(\cC^3)$ of $\cC^3/G$.
It should satisfy $\te_0<0$ and $\te_i>0$ for $1\le i\le5$.
A perfect matching $I\in\lA$ is \te-stable if and only if
one can reach any vertex of $Q$ from vertex $0\in Q_0$
by going only through the arrows of $Q\ms I$.
The \te-stable perfect matchings are listed in Table \ref{table:pm2}.
\begin{table}[hp]
\begin{tabular}{|c|c|c|}
\hline
$N$ &$I$& $\ub\hi_I$\\
\hline
4& $34,01,23,50,45,12$ & $e_1=(1,0,0)$\\
5& $34,14,30,50,52,12$ & $g^3=\frac16(3,0,3)$\\
7& $13,40,02,35,24,51$ & $e_2=(0,1,0)$\\
9& $13,40,23,50,24,51$ & $g^2=\frac16(2,4,0)$\\
10& $13,40,23,50,45,12$ & $g^4=\frac16(4,2,0)$\\
13& $41,40,30,50,52,51$ & $g=\frac16(1,2,3)$\\
17& $41,03,14,52,25,30$ & $e_3=(0,0,1)$\\
\hline
\end{tabular}
\bigskip
\caption{\te-stable perfect matchings}
\label{table:pm2}
\end{table}

We choose now such pairs of \te-stable perfect matchings
that their union is still \te-stable.
All pairs including $I_{13}$, except the pair
$\set{I_{13},I_4}$, satisfy this condition.
This allows us to reconstruct the toric diagram
of $\Hilb^G(\cC^3)$, see Figure \ref{toricd}.
This diagram coincides with a toric
diagram constructed by Nakamura \cite{Nakamura1}.
\begin{figure}[hp]
\scalebox{0.7}{\input{toricd1.pstex_t}}
\caption{Toric diagram of $\Hilb^{\cZ_6}(\cC^3)$}
\label{toricd}
\end{figure}

%% file: images/Z6-123.pstex_t
\begin{picture}(0,0)%
\includegraphics{Z6-123.pstex}%
\end{picture}%
\setlength{\unitlength}{4144sp}%
\begingroup\makeatletter\ifx\SetFigFont\undefined%
\gdef\SetFigFont#1#2#3#4#5{%
  \reset@font\fontsize{#1}{#2pt}%
  \fontfamily{#3}\fontseries{#4}\fontshape{#5}%
  \selectfont}%
\fi\endgroup%
\begin{picture}(6522,3839)(604,-3957)
\put(3961,-1141){\makebox(0,0)[lb]{\smash{{\SetFigFont{12}{14.4}{\familydefault}{\mddefault}{\updefault}{\color[rgb]{0,0,0}$4$}%
}}}}
\put(3061,-1141){\makebox(0,0)[lb]{\smash{{\SetFigFont{12}{14.4}{\familydefault}{\mddefault}{\updefault}{\color[rgb]{0,0,0}$1$}%
}}}}
\put(1261,-1141){\makebox(0,0)[lb]{\smash{{\SetFigFont{12}{14.4}{\familydefault}{\mddefault}{\updefault}{\color[rgb]{0,0,0}$1$}%
}}}}
\put(811,-3841){\makebox(0,0)[lb]{\smash{{\SetFigFont{12}{14.4}{\familydefault}{\mddefault}{\updefault}{\color[rgb]{0,0,0}$4$}%
}}}}
\put(811,-2041){\makebox(0,0)[lb]{\smash{{\SetFigFont{12}{14.4}{\familydefault}{\mddefault}{\updefault}{\color[rgb]{0,0,0}$3$}%
}}}}
\put(1711,-2041){\makebox(0,0)[lb]{\smash{{\SetFigFont{12}{14.4}{\familydefault}{\mddefault}{\updefault}{\color[rgb]{0,0,0}$0$}%
}}}}
\put(1261,-2941){\makebox(0,0)[lb]{\smash{{\SetFigFont{12}{14.4}{\familydefault}{\mddefault}{\updefault}{\color[rgb]{0,0,0}$2$}%
}}}}
\put(811,-3841){\makebox(0,0)[lb]{\smash{{\SetFigFont{12}{14.4}{\familydefault}{\mddefault}{\updefault}{\color[rgb]{0,0,0}$4$}%
}}}}
\put(1711,-3841){\makebox(0,0)[lb]{\smash{{\SetFigFont{12}{14.4}{\familydefault}{\mddefault}{\updefault}{\color[rgb]{0,0,0}$1$}%
}}}}
\put(2611,-3841){\makebox(0,0)[lb]{\smash{{\SetFigFont{12}{14.4}{\familydefault}{\mddefault}{\updefault}{\color[rgb]{0,0,0}$4$}%
}}}}
\put(2161,-2941){\makebox(0,0)[lb]{\smash{{\SetFigFont{12}{14.4}{\familydefault}{\mddefault}{\updefault}{\color[rgb]{0,0,0}$5$}%
}}}}
\put(3061,-2941){\makebox(0,0)[lb]{\smash{{\SetFigFont{12}{14.4}{\familydefault}{\mddefault}{\updefault}{\color[rgb]{0,0,0}$2$}%
}}}}
\put(2611,-2041){\makebox(0,0)[lb]{\smash{{\SetFigFont{12}{14.4}{\familydefault}{\mddefault}{\updefault}{\color[rgb]{0,0,0}$3$}%
}}}}
\put(7111,-2041){\makebox(0,0)[lb]{\smash{{\SetFigFont{12}{14.4}{\familydefault}{\mddefault}{\updefault}{\color[rgb]{0,0,0}$0$}%
}}}}
\put(4861,-2941){\makebox(0,0)[lb]{\smash{{\SetFigFont{12}{14.4}{\familydefault}{\mddefault}{\updefault}{\color[rgb]{0,0,0}$2$}%
}}}}
\put(7111,-3841){\makebox(0,0)[lb]{\smash{{\SetFigFont{12}{14.4}{\familydefault}{\mddefault}{\updefault}{\color[rgb]{0,0,0}$1$}%
}}}}
\put(6661,-2941){\makebox(0,0)[lb]{\smash{{\SetFigFont{12}{14.4}{\familydefault}{\mddefault}{\updefault}{\color[rgb]{0,0,0}$2$}%
}}}}
\put(6211,-2041){\makebox(0,0)[lb]{\smash{{\SetFigFont{12}{14.4}{\familydefault}{\mddefault}{\updefault}{\color[rgb]{0,0,0}$3$}%
}}}}
\put(3511,-3841){\makebox(0,0)[lb]{\smash{{\SetFigFont{12}{14.4}{\familydefault}{\mddefault}{\updefault}{\color[rgb]{0,0,0}$1$}%
}}}}
\put(6211,-3841){\makebox(0,0)[lb]{\smash{{\SetFigFont{12}{14.4}{\familydefault}{\mddefault}{\updefault}{\color[rgb]{0,0,0}$4$}%
}}}}
\put(5761,-2941){\makebox(0,0)[lb]{\smash{{\SetFigFont{12}{14.4}{\familydefault}{\mddefault}{\updefault}{\color[rgb]{0,0,0}$5$}%
}}}}
\put(5311,-2041){\makebox(0,0)[lb]{\smash{{\SetFigFont{12}{14.4}{\familydefault}{\mddefault}{\updefault}{\color[rgb]{0,0,0}$0$}%
}}}}
\put(3511,-2041){\makebox(0,0)[lb]{\smash{{\SetFigFont{12}{14.4}{\familydefault}{\mddefault}{\updefault}{\color[rgb]{0,0,0}$0$}%
}}}}
\put(4411,-2041){\makebox(0,0)[lb]{\smash{{\SetFigFont{12}{14.4}{\familydefault}{\mddefault}{\updefault}{\color[rgb]{0,0,0}$3$}%
}}}}
\put(5311,-3841){\makebox(0,0)[lb]{\smash{{\SetFigFont{12}{14.4}{\familydefault}{\mddefault}{\updefault}{\color[rgb]{0,0,0}$1$}%
}}}}
\put(3961,-2941){\makebox(0,0)[lb]{\smash{{\SetFigFont{12}{14.4}{\familydefault}{\mddefault}{\updefault}{\color[rgb]{0,0,0}$5$}%
}}}}
\put(4411,-3841){\makebox(0,0)[lb]{\smash{{\SetFigFont{12}{14.4}{\familydefault}{\mddefault}{\updefault}{\color[rgb]{0,0,0}$4$}%
}}}}
\put(6211,-241){\makebox(0,0)[lb]{\smash{{\SetFigFont{12}{14.4}{\familydefault}{\mddefault}{\updefault}{\color[rgb]{0,0,0}$2$}%
}}}}
\put(7111,-241){\makebox(0,0)[lb]{\smash{{\SetFigFont{12}{14.4}{\familydefault}{\mddefault}{\updefault}{\color[rgb]{0,0,0}$5$}%
}}}}
\put(6661,-1141){\makebox(0,0)[lb]{\smash{{\SetFigFont{12}{14.4}{\familydefault}{\mddefault}{\updefault}{\color[rgb]{0,0,0}$1$}%
}}}}
\put(4411,-241){\makebox(0,0)[lb]{\smash{{\SetFigFont{12}{14.4}{\familydefault}{\mddefault}{\updefault}{\color[rgb]{0,0,0}$2$}%
}}}}
\put(2611,-241){\makebox(0,0)[lb]{\smash{{\SetFigFont{12}{14.4}{\familydefault}{\mddefault}{\updefault}{\color[rgb]{0,0,0}$2$}%
}}}}
\put(811,-241){\makebox(0,0)[lb]{\smash{{\SetFigFont{12}{14.4}{\familydefault}{\mddefault}{\updefault}{\color[rgb]{0,0,0}$2$}%
}}}}
\put(1711,-241){\makebox(0,0)[lb]{\smash{{\SetFigFont{12}{14.4}{\familydefault}{\mddefault}{\updefault}{\color[rgb]{0,0,0}$5$}%
}}}}
\put(5311,-241){\makebox(0,0)[lb]{\smash{{\SetFigFont{12}{14.4}{\familydefault}{\mddefault}{\updefault}{\color[rgb]{0,0,0}$5$}%
}}}}
\put(3511,-241){\makebox(0,0)[lb]{\smash{{\SetFigFont{12}{14.4}{\familydefault}{\mddefault}{\updefault}{\color[rgb]{0,0,0}$5$}%
}}}}
\put(2161,-1141){\makebox(0,0)[lb]{\smash{{\SetFigFont{12}{14.4}{\familydefault}{\mddefault}{\updefault}{\color[rgb]{0,0,0}$4$}%
}}}}
\put(5761,-1141){\makebox(0,0)[lb]{\smash{{\SetFigFont{12}{14.4}{\familydefault}{\mddefault}{\updefault}{\color[rgb]{0,0,0}$4$}%
}}}}
\put(4861,-1141){\makebox(0,0)[lb]{\smash{{\SetFigFont{12}{14.4}{\familydefault}{\mddefault}{\updefault}{\color[rgb]{0,0,0}$1$}%
}}}}
\end{picture}%

%% file: images/Z6-123-mckay.pstex_t
\begin{picture}(0,0)%
\includegraphics{Z6-123-mckay.pstex}%
\end{picture}%
\setlength{\unitlength}{4144sp}%
\begingroup\makeatletter\ifx\SetFigFont\undefined%
\gdef\SetFigFont#1#2#3#4#5{%
  \reset@font\fontsize{#1}{#2pt}%
  \fontfamily{#3}\fontseries{#4}\fontshape{#5}%
  \selectfont}%
\fi\endgroup%
\begin{picture}(1920,2203)(2776,-1916)
\put(4681,-1276){\makebox(0,0)[lb]{\smash{{\SetFigFont{12}{14.4}{\familydefault}{\mddefault}{\updefault}{\color[rgb]{0,0,0}$2$}%
}}}}
\put(4681,-376){\makebox(0,0)[lb]{\smash{{\SetFigFont{12}{14.4}{\familydefault}{\mddefault}{\updefault}{\color[rgb]{0,0,0}$1$}%
}}}}
\put(2791,-421){\makebox(0,0)[lb]{\smash{{\SetFigFont{12}{14.4}{\familydefault}{\mddefault}{\updefault}{\color[rgb]{0,0,0}$5$}%
}}}}
\put(2791,-1321){\makebox(0,0)[lb]{\smash{{\SetFigFont{12}{14.4}{\familydefault}{\mddefault}{\updefault}{\color[rgb]{0,0,0}$4$}%
}}}}
\put(3871,164){\makebox(0,0)[lb]{\smash{{\SetFigFont{12}{14.4}{\familydefault}{\mddefault}{\updefault}{\color[rgb]{0,0,0}$0$}%
}}}}
\put(3871,-1861){\makebox(0,0)[lb]{\smash{{\SetFigFont{12}{14.4}{\familydefault}{\mddefault}{\updefault}{\color[rgb]{0,0,0}$3$}%
}}}}
\end{picture}%

%% file: images/toricd1.pstex_t
\begin{picture}(0,0)%
\includegraphics{toricd1.pstex}%
\end{picture}%
\setlength{\unitlength}{3947sp}%
\begingroup\makeatletter\ifx\SetFigFont\undefined%
\gdef\SetFigFont#1#2#3#4#5{%
  \reset@font\fontsize{#1}{#2pt}%
  \fontfamily{#3}\fontseries{#4}\fontshape{#5}%
  \selectfont}%
\fi\endgroup%
\begin{picture}(3173,3072)(1118,-4550)
\put(2851,-1636){\makebox(0,0)[lb]{\smash{{\SetFigFont{12}{14.4}{\familydefault}{\mddefault}{\updefault}{\color[rgb]{0,0,0}$e_2$}%
}}}}
\put(1276,-4486){\makebox(0,0)[lb]{\smash{{\SetFigFont{12}{14.4}{\familydefault}{\mddefault}{\updefault}{\color[rgb]{0,0,0}$e_3$}%
}}}}
\put(2776,-4486){\makebox(0,0)[lb]{\smash{{\SetFigFont{12}{14.4}{\familydefault}{\mddefault}{\updefault}{\color[rgb]{0,0,0}$g^3$}%
}}}}
\put(4276,-4486){\makebox(0,0)[lb]{\smash{{\SetFigFont{12}{14.4}{\familydefault}{\mddefault}{\updefault}{\color[rgb]{0,0,0}$e_1$}%
}}}}
\put(2176,-3586){\makebox(0,0)[lb]{\smash{{\SetFigFont{12}{14.4}{\familydefault}{\mddefault}{\updefault}{\color[rgb]{0,0,0}$g$}%
}}}}
\put(3826,-3436){\makebox(0,0)[lb]{\smash{{\SetFigFont{12}{14.4}{\familydefault}{\mddefault}{\updefault}{\color[rgb]{0,0,0}$g^4$}%
}}}}
\put(3376,-2536){\makebox(0,0)[lb]{\smash{{\SetFigFont{12}{14.4}{\familydefault}{\mddefault}{\updefault}{\color[rgb]{0,0,0}$g^2$}%
}}}}
\end{picture}%

%% file: a12.bbl
\providecommand{\bysame}{\leavevmode\hbox to3em{\hrulefill}\thinspace}
\providecommand{\href}[2]{#2}
\begin{thebibliography}{10}

\bibitem{Aspinwall}
Paul~S. Aspinwall, \emph{{D}-branes on toric {C}alabi-{Y}au varieties},
  \eprint{arxiv}{0806.2612}.

\bibitem{AtiyahMacdonald}
M.~F. Atiyah and I.~G. Macdonald, \emph{Introduction to commutative algebra},
  Addison-Wesley Publishing Co., Reading, Mass.-London-Don Mills, Ont., 1969.

\bibitem{Bass1}
Hyman Bass, \emph{On the ubiquity of {G}orenstein rings}, Math. Z. \textbf{82}
  (1963), 8--28.

\bibitem{Bock2}
Raf Bocklandt, \emph{Calabi {Y}au algebras and weighted quiver polyhedra},
  \eprint{arxiv}{0905.0232}.

\bibitem{Bock1}
\bysame, \emph{Graded {C}alabi-{Y}au algebras of dimension $3$}, J. Pure Appl.
  Algebra \textbf{212} (2008), no.~1, 14--32, \eprint{arxiv}{math.RA/0603558}.

\bibitem{BKR}
Tom Bridgeland, Alastair King, and Miles Reid, \emph{The {M}c{K}ay
  correspondence as an equivalence of derived categories}, J. Amer. Math. Soc.
  \textbf{14} (2001), no.~3, 535--554, \eprint{arxiv}{math.AG/9908027}.

\bibitem{Broomhead1}
Nathan Broomhead, \emph{Dimer models and {Calabi-Yau} algebras},
  \eprint{arxiv}{0901.4662}, {PhD Thesis}.

\bibitem{HerzogBruns}
Winfried Bruns and J{\"u}rgen Herzog, \emph{Cohen-{M}acaulay rings}, Cambridge
  Studies in Advanced Mathematics, vol.~39, Cambridge University Press,
  Cambridge, 1993.

\bibitem{Cartan1}
Henri Cartan and Samuel Eilenberg, \emph{Homological algebra}, Princeton
  University Press, Princeton, N. J., 1956.

\bibitem{Craw2}
Alastair Craw, Diane Maclagan, and Rekha~R. Thomas, \emph{Moduli of {M}c{K}ay
  quiver representations. {I}. {T}he coherent component}, Proc. Lond. Math.
  Soc. (3) \textbf{95} (2007), no.~1, 179--198,
  \eprint{arxiv}{math.AG/0505115}.

\bibitem{Craw3}
\bysame, \emph{Moduli of {M}c{K}ay quiver representations. {II}. {G}r\"obner
  basis techniques}, J. Algebra \textbf{316} (2007), no.~2, 514--535,
  \eprint{arxiv}{math.AG/0611840}.

\bibitem{Davison1}
Ben Davison, \emph{Consistency conditions for brane tilings},
  \eprint{arxiv}{0812.4185}.

\bibitem{Ginz1}
Victor Ginzburg, \emph{Calabi-{Y}au algebras}, \eprint{arxiv}{math.AG/0612139}.

\bibitem{HHV}
Amihay Hanany, Christopher~P. Herzog, and David Vegh, \emph{Brane tilings and
  exceptional collections}, J. High Energy Phys. (2006), no.~7, 001, 44 pp.
  (electronic), \eprint{arxiv}{hep-th/0602041v2}.

\bibitem{Ishii2}
Akira Ishii and Kazushi Ueda, \emph{Dimer models and the special {McKay}
  correspondence}, \eprint{arxiv}{0905.0059}.

\bibitem{Ishii1}
\bysame, \emph{On moduli spaces of quiver representations associated with brane
  tilings}, \eprint{arxiv}{0710.1898}.

\bibitem{Ito-Nakajima1}
Yukari Ito and Hiraku Nakajima, \emph{Mc{K}ay correspondence and {H}ilbert
  schemes in dimension three}, Topology \textbf{39} (2000), no.~6, 1155--1191,
  \eprint{arxiv}{math.AG/9803120}.

\bibitem{King1}
A.~D. King, \emph{Moduli of representations of finite-dimensional algebras},
  Quart. J. Math. Oxford Ser. (2) \textbf{45} (1994), no.~180, 515--530.

\bibitem{MR2}
Sergey Mozgovoy and Markus Reineke, \emph{On the non-commutative
  {D}onaldson-{T}homas invariants arising from brane tilings},
  \eprint{arxiv}{0809.0117}.

\bibitem{Nakamura1}
Iku Nakamura, \emph{Hilbert schemes of abelian group orbits}, J. Algebraic
  Geom. \textbf{10} (2001), no.~4, 757--779.

\bibitem{Oda1}
Tadao Oda, \emph{Convex bodies and algebraic geometry}, Ergebnisse der
  Mathematik und ihrer Grenzgebiete (3) [Results in Mathematics and Related
  Areas (3)], vol.~15, Springer-Verlag, Berlin, 1988, An introduction to the
  theory of toric varieties, Translated from the Japanese.

\bibitem{Ooishi}
Akira Ooishi, \emph{Matlis duality and the width of a module}, Hiroshima Math.
  J. \textbf{6} (1976), no.~3, 573--587.

\bibitem{sardoinfirri-1996}
Alexander~V. Sardo-Infirri, \emph{Resolutions of orbifold singularities and
  flows on the {McKay} quiver}, \eprint{arxiv}{alg-geom/9610005}, Preprint.

\bibitem{Schrijver1}
Alexander Schrijver, \emph{Theory of linear and integer programming},
  Wiley-Interscience Series in Discrete Mathematics, John Wiley \& Sons Ltd.,
  Chichester, 1986, A Wiley-Interscience Publication.

\bibitem{Thaddeus1}
Michael Thaddeus, \emph{Toric quotients and flips}, Topology, geometry and
  field theory, World Sci. Publ., River Edge, NJ, 1994, pp.~193--213.

\bibitem{Bergh5}
Michel van~den Bergh, \emph{Existence theorems for dualizing complexes over
  non-commutative graded and filtered rings}, J. Algebra \textbf{195} (1997),
  no.~2, 662--679.

\bibitem{Bergh1}
\bysame, \emph{Non-commutative crepant resolutions}, The legacy of {N}iels
  {H}enrik {A}bel, Springer, Berlin, 2004, \eprint{arxiv}{math/0211064},
  pp.~749--770.

\bibitem{Yek1}
Amnon Yekutieli, \emph{Dualizing complexes, {M}orita equivalence and the
  derived {P}icard group of a ring}, J. London Math. Soc. (2) \textbf{60}
  (1999), no.~3, 723--746, \eprint{arxiv}{math.RA/9810134}.

\end{thebibliography}
